# ON ASYMPTOTICS OF EIGENVECTORS OF LARGE SAMPLE COVARIANCE MATRIX

By Z. D. Bai,[1,2] B. Q. Miao[3] and G. M. Pan[2,3]

*Northeast Normal University, National University of Singapore and University of Science and Technology of China*

Let $\{X_{ij}\}$, $i,j=\dots$, be a double array of i.i.d. complex random variables with $EX_{11}=0, E|X_{11}|^2=1$ and $E|X_{11}|^4<\infty$, and let $A_n=\frac{1}{N}T_n^{1/2}X_nX_n^*T_n^{1/2}$, where $T_n^{1/2}$ is the square root of a nonnegative definite matrix $T_n$ and $X_n$ is the $n\times N$ matrix of the upper-left corner of the double array. The matrix $A_n$ can be considered as a sample covariance matrix of an i.i.d. sample from a population with mean zero and covariance matrix $T_n$, or as a multivariate $F$ matrix if $T_n$ is the inverse of another sample covariance matrix. To investigate the limiting behavior of the eigenvectors of $A_n$, a new form of empirical spectral distribution is defined with weights defined by eigenvectors and it is then shown that this has the same limiting spectral distribution as the empirical spectral distribution defined by equal weights. Moreover, if $\{X_{ij}\}$ and $T_n$ are either real or complex and some additional moment assumptions are made then linear spectral statistics defined by the eigenvectors of $A_n$ are proved to have Gaussian limits, which suggests that the eigenvector matrix of $A_n$ is nearly Haar distributed when $T_n$ is a multiple of the identity matrix, an easy consequence for a Wishart matrix.

**1. Introduction.** Let $X_n=(X_{ij})$ be an $n\times N$ matrix of i.i.d. complex random variables and let $T_n$ be an $n\times n$ nonnegative definite Hermitian

Received June 2005; revised March 2006.
[1]Supported by NSFC Grant 10571020.
[2]Supported in part by NUS Grant R-155-000-056-112.
[3]Supported by Grants 10471135 and 10571001 from the National Natural Science Foundation of China.
*AMS 2000 subject classifications.* Primary 15A52, 60F15, 62E20; secondary 60F17, 62H99.
*Key words and phrases.* Asymptotic distribution, central limit theorems, CDMA, eigenvectors and eigenvalues, empirical spectral distribution function, Haar distribution, MIMO, random matrix theory, sample covariance matrix, SIR, Stieltjes transform, strong convergence.

matrix with a square root $T_n^{1/2}$. In this paper, we shall consider the matrix $A_n = \frac{1}{N} T_n^{1/2} X_n X_n^* T_n^{1/2}$. If $T_n$ is nonrandom, then $A_n$ can be considered as a sample covariance matrix of a sample drawn from a population with the same distribution as $T_n^{1/2}\mathbf{X}_{\cdot,1}$, where $\mathbf{X}_{\cdot,1} = (X_{11}, \ldots, X_{n1})'$. If $T_n$ is an inverse of another sample covariance matrix, then the multivariate $F$ matrix can be considered as a special case of the matrix $A_n$.

In this paper, we consider the case where both dimension $n$ and sample size $N$ are large. Bai and Silverstein [7] gave an example demonstrating the considerable difference between the case where $n$ is fixed and that where $n$ increases with $N$ proportionally. When $T_n = I$, $A_n$ reduces to the usual sample covariance matrix of $N$ $n$-dimensional random vectors with mean 0 and covariance matrix $I$. An important statistic in multivariate analysis is

$$W_n = \ln(\det A_n) = \sum_{j=1}^{N} \ln(\lambda_j),$$

where $\lambda_j, j = 1, \ldots, n$, are the eigenvalues of $A_n$. When $n$ is fixed, by taking a Taylor expression of $\ln(1+x)$, one can easily prove that

$$\sqrt{\frac{N}{n}} W_n \xrightarrow{D} N(0, EX_{11}^4 - 1).$$

It appears that when $n$ is fixed, this distribution can be used to test the hypothesis of variance homogeneity. However, it is not the case when $n$ increases as $[cN]$ (the integer part of $cN$) with $c \in (0,1)$. Using results of the limiting spectral distribution of $A_n$ (see [12] or [1]), one can show that with probability one that

$$\frac{1}{n} W_n \to \frac{c-1}{c} \ln(1-c) - 1 \equiv d(c) < 0,$$

which implies that

$$\sqrt{\frac{N}{n}} W_n \sim d(c)\sqrt{Nn} \to -\infty.$$

More precisely, the distribution of $W_n$ shaft to left quickly when $n$ increases as $n \sim cN$. Figure 1 gives the kernel density estimates using 1000 realizations of $W_n$ for $N = 20, 100, 200$ and 500 with $n = 0.2N$. Figures 2 and 3 give the kernel density estimates of $\sqrt{\frac{N}{n}} W_n$ for the cases $n = 5$ and $n = 10$ with $N = 50$. These figures clearly show that the distribution of $W_n$ cannot be approximated by a centered normal distribution even if the ratio $c$ is as small as 0.1.

This phenomenon motivates the development of the theory of spectral analysis of large-dimensional random matrices which is simply called random matrix theory (RMT). In this theory, for a square matrix $A$ of real eigenvalues, its empirical spectral distribution (ESD) $F^A$ is defined as the empirical distribution generated by its eigenvalues. The limiting properties of the ESD of sample covariance matrices have been intensively investigated in the literature and the reader is referred to [[1, 5, 6, 7, 11, 12, 13, 17, 20, 21, 23]].

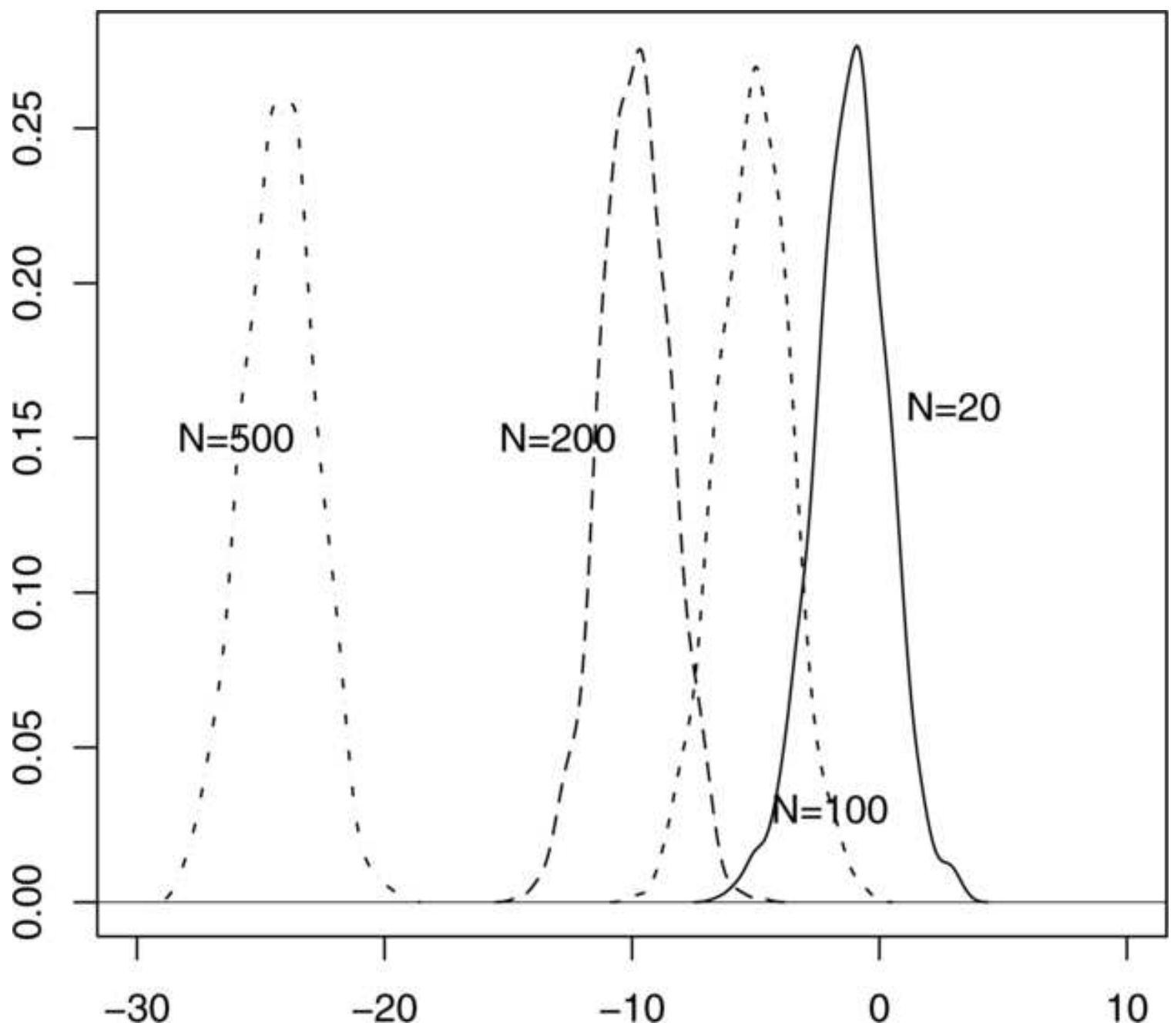

FIG. 1. *Density of $W_n$ under different sample sizes, $c = 0.2$.*

An important mathematical tool in RMT is the *Stieltjes transform*, which is defined by

$$m_G(z) = \int \frac{1}{\lambda - z} \, dG(\lambda), \qquad z \in \mathbb{C}^+ \equiv \{z \in \mathbb{C}, \ \Im z > 0\},$$

for any distribution function $G(x)$. It is well known that $G_n \xrightarrow{v} G$ if and only if $m_{G_n}(z) \to m_G(z)$, for all $z \in \mathbb{C}^+$.

In [13], it is assumed that:

(i) for all $n, i, j$, $X_{ij}$ are independently and identically distributed with $EX_{ij} = 0$ and $E|X_{ij}|^2 = 1$;

(ii) $F^{T_n} \xrightarrow{\mathcal{D}} H$, a proper distribution function;

(iii) $\frac{n}{N} \to c > 0$ as $n \to \infty$.

It is then proved that, with probability 1, $F^{A_n}$ converges to a nonrandom distribution function $F^{c,H}$ whose Stieltjes transform $m(z) = m_{F^{c,H}}(z)$, for

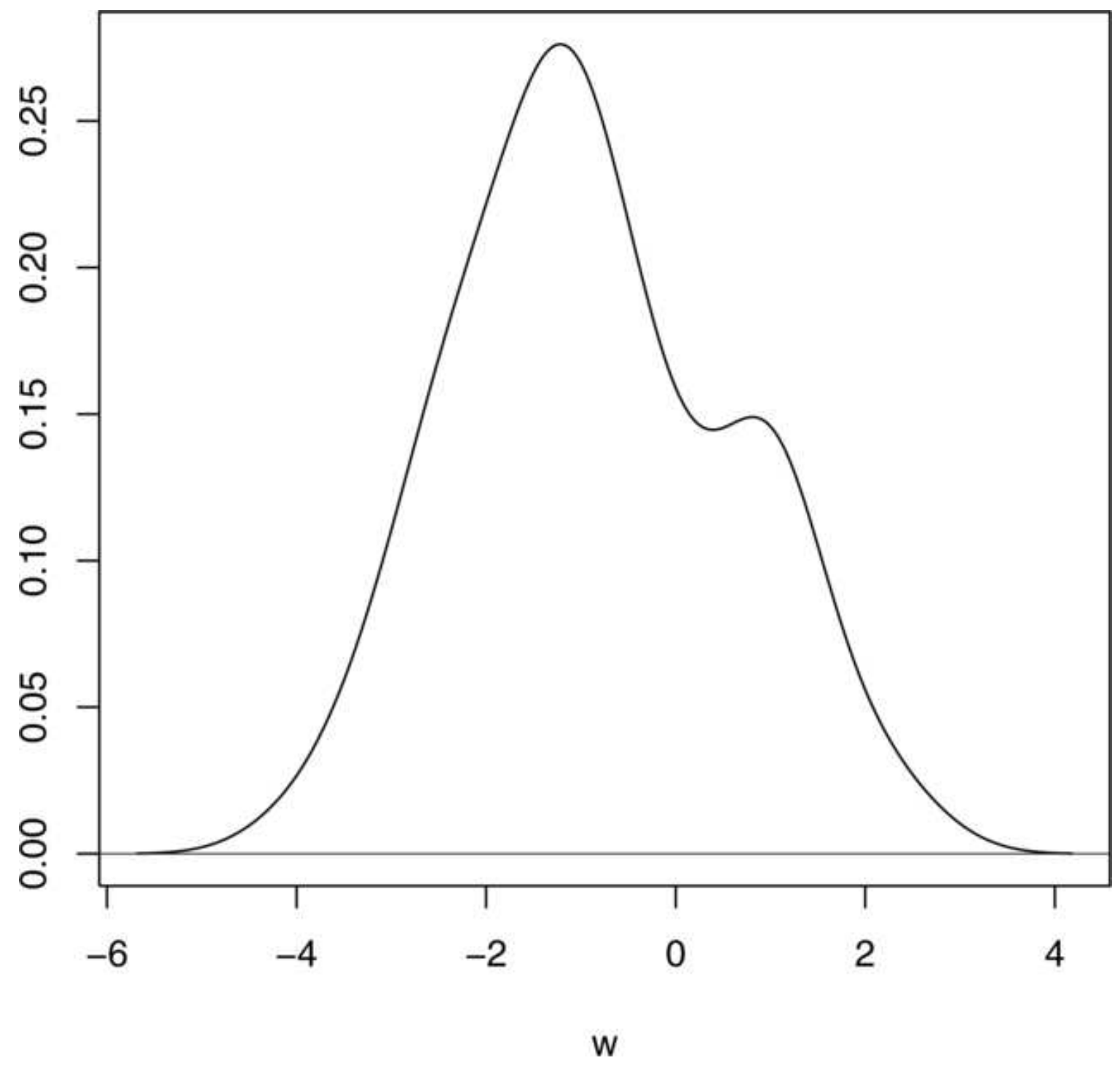

FIG. 2. *Probability density of* $\sqrt{\frac{N}{n}}W_n$ $(n=5, N=50)$.

each $z \in \mathbb{C}^+$, is the unique solution in $\mathbb{C}^+$ of the equation

$$m(z) = \int \frac{1}{t(1-c-czm)-z}\, dH(t). \tag{1.1}$$

Let $\underline{A}_n = \frac{1}{N} X_n^* T_n X_N$. The spectrum of $A_n$ differs from that of $\underline{A}_n$ only by $|n-N|$ zero eigenvalues. Hence, we have

$$F^{\underline{A}_n} = \left(1 - \frac{n}{N}\right) I_{[0,\infty)} + \frac{n}{N} F^{A_n}.$$

It then follows that

$$\underline{m}_n(z) = m_{F^{\underline{A}_n}}(z) = -\frac{1-n/N}{z} + \frac{n}{N} m_{F^{A_n}}(z) \tag{1.2}$$

and, correspondingly for their limits,

$$\underline{m}(z) = m_{\underline{F}^{c,H}}(z) = -\frac{1-c}{z} + cm(z), \tag{1.3}$$

where $\underline{F}^{c,H}$ is the limiting empirical distribution function of $\underline{A}_n$.

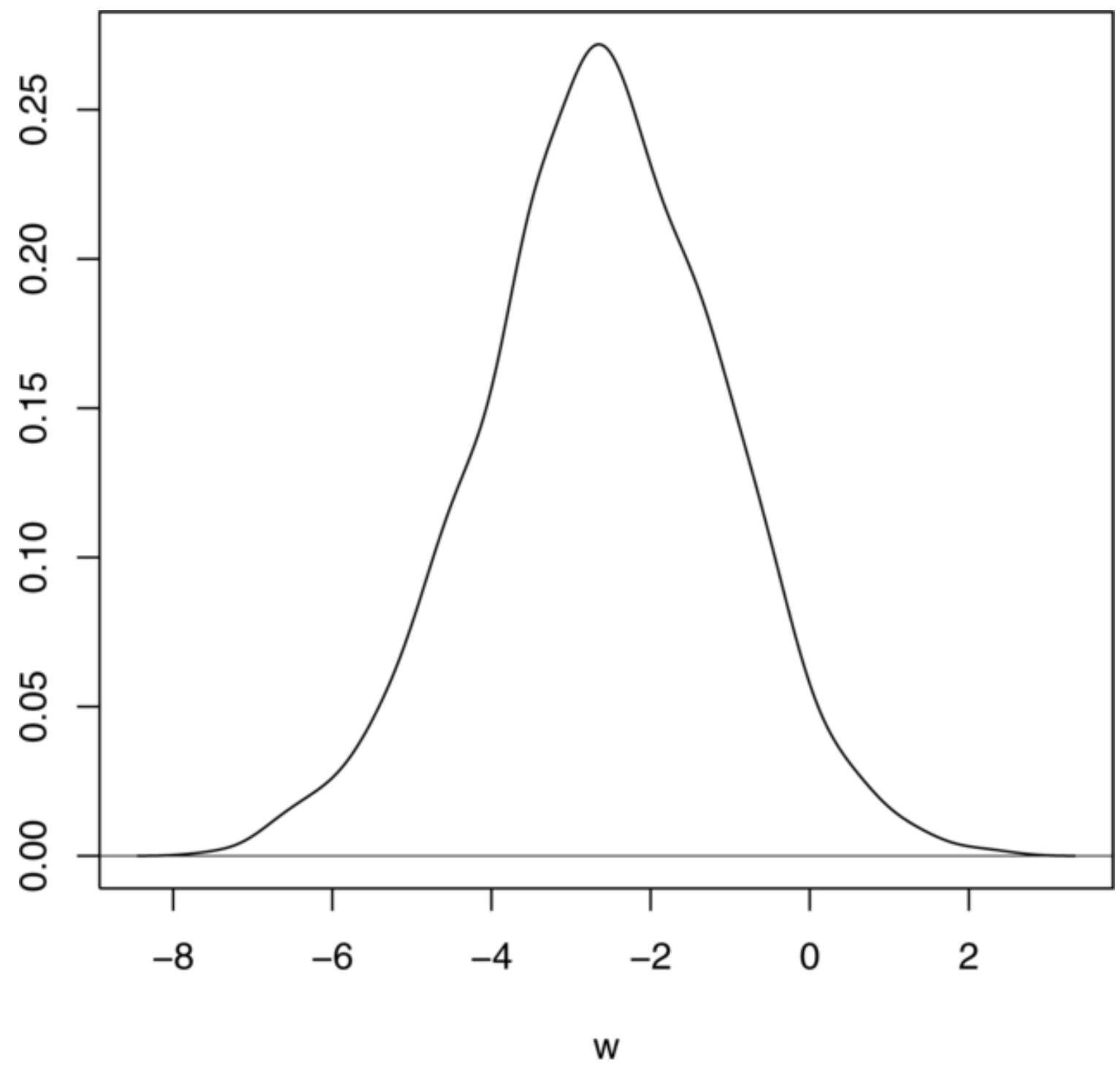

FIG. 3. *Probability density of* $\sqrt{\frac{N}{n}}W_n$ $(n = 10, N = 50)$.

Using this notation, equation (1.1) can be converted to equation (1.4) for $\underline{m}(z)$. That is, $\underline{m}(z)$ is the unique solution in $\mathbb{C}^+$ of the equation

$$\underline{m} = -\left(z - c\int \frac{t\,dH(t)}{1 + t\underline{m}}\right)^{-1}. \tag{1.4}$$

From this equation, the inverse function has the explicit form

$$z = -\frac{1}{\underline{m}} + c\int \frac{t\,dH(t)}{1 + t\underline{m}}. \tag{1.5}$$

The limiting properties of the eigenvalues of $A_n$ have been intensively investigated in the literature. Among others, we shall now briefly mention some remarkable ones. Yin, Bai and Krishnaiah [22] established the limiting value of the largest eigenvalue, while Bai and Yin [2] employed a unified approach to obtain the limiting value for the smallest and largest eigenvalues of $A_n$ when $T_N = I$. A breakthrough on the convergence rate of the ESD of a sample covariance matrix was made in [3] and [4]. In [5], it is shown that, with probability 1, no eigenvalues of $A_n$ appear in any interval $[a, b]$ which is contained in an open interval outside the supports of $F^{c_n,H_n}$ for all large $N$ under the condition of finite 4th moment (here, $c_n = n/N$ and $H_n$ is the ESD of $T_n$).

However, relatively less work has been done on the limiting behavior of eigenvectors of $A_n$. Some results on this aspect can be found in [[14, 15, 16]]. That more attention has been paid to the ESD of the sample covariance matrix may be due to the origins of RMT, which lie with quantum mechanics (QM), where the eigenvalues of large-dimensional random matrices are used to describe energy levels of particles. With the application of RMT to many other areas, such as statistics, wireless communications,—for example, the CDMA (code division multiple access) systems and MIMO (multiple input multiple output) systems, finance and economics, and so on, the importance of the limiting behavior of eigenvectors has been gradually recognized. For example, in signal processing, for signals received by linearly spaced sensors, the estimates of the directions of arrivals (DOA) are based on the noise eigenspace. In principal component analysis or factor analysis, the directions of the principal components are the eigenvectors corresponding to the largest eigenvalues.

Now, let us consider another example of the application of eigenvectors of a large covariance matrix $A_n$ which is important in wireless communications. Consider transmission methods in wireless systems. In a direct sequence CDMA system, the discrete-time model for a synchronous systems is formulated as

$$\mathbf{r} = \sum_{k=1}^{K} b_k \mathbf{s}_k + \mathbf{w},$$

where $b_k (\in \mathbb{C})$ and $\mathbf{s}_k (\in \mathbb{C}^N)$ are the transmitted data symbols and signature sequence of the user $k$, respectively, and $\mathbf{w}$ is an $N$-dimensional background Gaussian noise of i.i.d. variables with mean zero and variance $\sigma^2$. The goal is to demodulate the transmitted $b_k$ for each user. In this case, the performance measure is defined as the signal-to-interference ratio (SIR) of the estimates. In a large network, since the number of users is very large, it is reasonable to assume that $K$ is proportional to $N$. That is, one can assume that their ratio remains constant when both $K$ and $N$ tend to infinity. Thus, it is feasible to apply the theory of large-dimensional random matrices to wireless communications and, indeed, there has already accumulated a fruitful literature in this direction (see, e.g., [18] and [19], among others). Eldar and Chen [10] derived an expression of SIR for the decorrelator receiver in terms of eigenvectors and eigenvalues of random matrices and then analyzed the asymptotics of the SIR (see [10] for details).

Our research is motivated by the fact that the matrix of eigenvectors (eigenmatrix for short) of the Wishart matrix has the Haar distribution, that is, the uniform distribution over the group of unitary matrices (or orthogonal matrices in the real case). It is conceivable that the eigenmatrix of a large sample covariance matrix should be "asymptotically Haar distributed." However, we are facing a problem on how to formulate the terminology "asymptotically Haar distributed" because the dimensions of the eigenmatrices are increasing. In this paper, we shall adopt the method of Silverstein [14, 15]. If $U$ has a Haar measure over the orthogonal matrices, then for any unit vector $\mathbf{x} \in \mathbb{R}^n$, $\mathbf{y} = U\mathbf{x} = (y_1, \ldots, y_n)'$ has a uniform distribution over the unit sphere $S_n = \{\mathbf{x} \in \mathbb{R}^n;\ \|\mathbf{x}\| = 1\}$. If $\mathbf{z} = (z_1, \ldots, z_n)' \sim N(0, I_n)$, then $\mathbf{y}$ has the same distribution as $\mathbf{z}/\|\mathbf{z}\|$.

Now, define a stochastic process $Y_n(t)$ in the space $D(0,1)$ by

$$\begin{aligned} Y_n(t) &= \sqrt{\frac{n}{2}} \sum_{i=1}^{[nt]} \left( |y_i|^2 - \frac{1}{n} \right) \\ &\stackrel{d}{=} \sqrt{\frac{n}{2}} \frac{1}{\|\mathbf{z}\|^2} \sum_{i=1}^{[nt]} \left( |z_i|^2 - \frac{\|\mathbf{z}\|^2}{n} \right), \end{aligned} \tag{1.6}$$

where $[a]$ denotes the greatest integer $\le a$. From the second equality, it is easy to see that $Y_n(t)$ converges to a Brownian bridge (BB) $B(t)$ when $n$ converges to infinity. Thus we are interested in whether the same is true for general sample covariance matrices.

Let $U_n \Lambda_n U_n^*$ denote the spectral decomposition of $A_n$, where $\Lambda_n = \operatorname{diag}(\lambda_1, \lambda_2, \ldots, \lambda_n)$ and $U_n = (u_{ij})$ is a unitary matrix consisting of the orthonormal eigenvectors of $A_n$. Assume that $\mathbf{x}_n \in \mathbb{C}^n$, $\|\mathbf{x}_n\| = 1$, is an arbitrary nonrandom unit vector and that $\mathbf{y} = (y_1, y_2, \ldots, y_n)^* = U_n^* \mathbf{x}_n$. We define a stochastic process by way of (1.6). If $U_n$ is "asymptotically Haar distributed," then $\mathbf{y}$ should be "asymptotically uniformly distributed" over the unit sphere and $Y_n(t)$ should tend to a BB. Our main goal is to examine the limiting properties of the vector $\mathbf{y}$ through the stochastic process $Y_n(t)$.

For ease of application of RMT, we make a time transformation in $Y_n(t)$,

$$X_n(x) = Y_n(F^{A_n}(x)),$$

where $F^{A_n}$ is the ESD of the matrix $A_n$. If the distribution of $U_n$ is close to Haar, then $X_n(x)$ should approximate $B(F^{c,H}(x))$, where $F^{c,H}$ is the limiting spectral distribution of $A_n$.

We define a new empirical distribution function based on eigenvectors and eigenvalues:

$$F_1^{A_n}(x) = \sum_{i=1}^{n} |y_i|^2 I(\lambda_i \le x). \tag{1.7}$$

Recall that the ESD of $A_n$ is

$$F^{A_n}(x) = \frac{1}{n} \sum_{i=1}^{n} I(\lambda_i \le x).$$

It then follows that

$$X_n(x) = \sqrt{\frac{n}{2}}(F_1^{A_n}(x) - F^{A_n}(x)).$$

The investigation of $Y_n(t)$ is then converted to one concerning the difference $F_1^{A_n}(x) - F^{A_n}(x)$ of the two empirical distributions.

It is obvious that $F_1^{A_n}(x)$ is a random probability distribution function and that its Stieltjes transform is given by

$$m_{F_1^{A_n}}(z) = x_n^*(A_n - zI)^{-1}x_n. \tag{1.8}$$

As we have seen, the difference between $F_1^{A_n}(x)$ and $F^{A_n}(x)$ is only in their different weights on the eigenvalues of $A_n$. However, it will be proven that although these two empirical distributions have different weights, they have the same limit; this is included in Theorem 1.1 below.

To investigate the convergence of $X_n(x)$, we consider its linear functional, which is defined as

$$\hat{X}_n(g) = \int g(x)\,dX_n(x),$$

where $g$ is a bounded continuous function. It turns out that

$$\begin{aligned}\hat{X}_n(g) &= \sqrt{\frac{n}{2}}\left[\sum_{j=1}^n |y_j^2| g(\lambda_j) - \frac{1}{n}\sum_{j=1}^n g(\lambda_j)\right]\\ &= \sqrt{\frac{n}{2}}\left[\int g(x)\,dF_1^{A_n}(x) - \int g(x)\,dF^{A_n}(x)\right].\end{aligned}$$

Proving the convergence of $\hat{X}_n(g)$ under general conditions is difficult. Following an idea of [7], we shall prove the central limit theorem (CLT) for those $g$ which are analytic over the support of the limiting spectral distribution of $A_n$. To this end, let

$$G_n(x) = \sqrt{N}(F_1^{A_n}(x) - F^{c_n,H_n}(x)),$$

where $c_n = \frac{n}{N}$ and where $F^{c_n,H_n}(x)$ denotes the limiting distribution by substituting $c_n$ for $c$ and $H_n$ for $H$ in $F^{c,H}$.

The main results of this paper are formulated in the following three theorems.

THEOREM 1. *Suppose that:*

(1) *for each $n$, $X_{ij} = X_{ij}^n, i,j = 1,2,\ldots,$ are i.i.d. complex random variables with $EX_{11} = 0, E|X_{11}|^2 = 1$ and $E|X_{11}|^4 < \infty$;*

(2) $\mathbf{x}_n \in \mathbb{C}_1^n = \{\mathbf{x} \in \mathbb{C}^n, \|\mathbf{x}\| = 1\}$ *and* $\lim_{n\to\infty} \frac{n}{N} = c \in (0,\infty)$;

(3) $T_n$ *is nonrandom Hermitian nonnegative definite with its spectral norm bounded in* $n$, *with* $H_n = F^{T_n} \xrightarrow{\mathcal{D}} H$ *a proper distribution function and with* $\mathbf{x}_n^*(T_n - zI)^{-1}\mathbf{x}_n \to m_{F^H}(z)$, *where* $m_{F^H}(z)$ *denotes the Stieltjes transform of* $H(t)$. *It then follows that*

$$F_1^{A_n}(x) \to F^{c,H}(x) \qquad a.s.$$

REMARK 1. The condition $\mathbf{x}_n^*(T_n - zI)^{-1}\mathbf{x}_n \to m_{F^H}(z)$ is critical for our Theorem 1 as well as for the main theorems which we give later. At first, we indicate that if $T_n = bI$ for some positive constant $b$ or, more generally, $\lambda_{\max}(T_n) - \lambda_{\min}(T_n) \to 0$, then the condition $\mathbf{x}_n^*(T_n - zI)^{-1}\mathbf{x}_n \to m_{F^H}(z)$ holds uniformly for all $\mathbf{x}_n \in \mathbb{C}_1^n$.

We also note that this condition does not require $T_n$ to be a multiple of an identity. As an application of this remark, one sees that the eigenmatrix of a sample covariance matrix transforms $\mathbf{x}_n$ to a unit vector whose entries' absolute values are close to $1/\sqrt{N}$. Consequently, the condition $\mathbf{x}_n^*(T_n - zI)^{-1}\mathbf{x}_n \to m_{F^H}(z)$ holds when $T_n$ is the inverse of another sample covariance matrix which is independent of $X_n$. Therefore, the multivariate $F$ matrix satisfies Theorem 1.

In general, the condition may not hold for all $\mathbf{x}_n \in \mathbb{C}_1^n$. However, there always exist some $\mathbf{x}_n \in \mathbb{C}_1^n$ such that this condition holds, say $\mathbf{x}_n = (\mathbf{u}_1 + \cdots + \mathbf{u}_n)/\sqrt{n}$, where $\mathbf{u}_1, \ldots, \mathbf{u}_n$ are the orthonormal eigenvectors of the spectral decomposition of $T_n$.

Applying Theorem 1, we get the following interesting results.

COROLLARY 1. *Let* $(A_n^m)_{ii}, m = 1, 2, \ldots,$ *denote the* $i$*th diagonal elements of matrices* $A_n^m$. *Under the conditions of Theorem* 1 *for* $\mathbf{x}_n = \mathbf{e}_{ni}$, *it follows that for any fixed* $m$,

$$\lim_{n\to\infty}\left|(A_n^m)_{ii} - \int x^m \, dF^{c,H}(x)\right| \to 0 \qquad a.s.,$$

*where* $\mathbf{e}_{ni}$ *is the* $n$*-vector with* $i$*th element* 1 *and all others* 0.

REMARK 2. If $T_n = bI$ for some positive constant $b$ or, more generally, $\lambda_{\max}(T_n) - \lambda_{\min}(T_n) \to 0$, then there is a better result, that is,

$$\lim_{n\to\infty}\max_i\left|(A_n^m)_{ii} - \int x^m \, dF^{c,H}(x)\right| \to 0 \qquad \text{a.s.} \tag{1.9}$$

[The corollary follows easily from Theorem 1. The uniform convergence of (1.9) follows from the uniform convergence of condition (3) of Theorem 1 and by careful checking of the proof of Theorem 1.]

More generally, we have the following:

COROLLARY 2. *If $f(x)$ is a bounded function and the assumptions of Theorem 1 are satisfied, then*

$$\sum_{j=1}^{n} |y_j^2| f(\lambda_j) - \frac{1}{n}\sum_{j=1}^{n} f(\lambda_j) \to 0 \qquad a.s.$$

REMARK 3. The proof of the above corollaries are immediate. Applying Corollary 2, Theorem 1 of [10] can be extended to a more general case without difficulty.

THEOREM 2. *In addition to the conditions of Theorem 1, we further assume that:*

(4) *$g_1, \dots, g_k$ are defined and analytic on an open region $\mathcal{D}$ of the complex plane which contains the real interval*

$$\left[\liminf_n \lambda_{\min}^{T_n} I_{(0,1)}(c)(1-\sqrt{c})^2, \limsup_n \lambda_{\max}^{T_n}(1+\sqrt{c})^2\right] \tag{1.10}$$

*and*

(5)

$$\sup_z \sqrt{N}\left| x_n^*(\underline{m}_{F^{c_n,H_n}}(z)T_n + I)^{-1}x_n - \int \frac{1}{\underline{m}_{F^{c_n,H_n}}(z)t+1}\,dH_n(t)\right| \to 0$$

*as $n \to \infty$.*

*Then the following conclusions hold:*

(a) *The random vectors*

$$\left(\int g_1(x)\,dG_n(x), \dots, \int g_k(x)\,dG_n(x)\right) \tag{1.11}$$

*form a tight sequence.*

(b) *If $X_{11}$ and $T_n$ are real and $EX_{11}^4 = 3$, then the above random vector converges weakly to a Gaussian vector $X_{g_1}, \dots, X_{g_k}$ with mean zero and covariance function*

$$\begin{aligned} &\operatorname{Cov}(X_{g_1}, X_{g_2}) \\ &\quad = -\frac{1}{2\pi^2}\int_{\mathcal{C}_1}\int_{\mathcal{C}_2} g_1(z_1)g_2(z_2) \\ &\qquad\qquad \times \frac{(z_2\underline{m}(z_2) - z_1\underline{m}(z_1))^2}{c^2 z_1 z_2 (z_2 - z_1)(\underline{m}(z_2) - \underline{m}(z_1))}\,dz_1\,dz_2. \end{aligned} \tag{1.12}$$

*The contours $\mathcal{C}_1$ and $\mathcal{C}_2$ in the above equation are disjoint, are both contained in the analytic region for the functions $(g_1, \dots, g_k)$ and both enclose the support of $F^{c_n,H_n}$ for all large $n$.*

(c) *If* $X_{11}$ *is complex with* $EX_{11}^2 = 0$ *and* $E|X_{11}|^4 = 2$, *then the conclusions* (a) *and* (b) *still hold, but the covariance function reduces to half of the quantity given in* (1.12).

REMARK 4. If $T_n = bI$ for some positive constant $b$ or, more generally, $\sqrt{n}(\lambda_{\max}(T_n) - \lambda_{\min}(T_n)) \to 0$, then condition (5) holds uniformly for all $\mathbf{x}_n \in \mathbb{C}_1^n$.

REMARK 5. Indeed, we can also establish the central limit theorem for $\hat{X}_n(g)$ according to Theorem 1.1 of [7] and Theorem 2. Beside Theorem 2, which holds for more general functions $g(x)$, the following theorem reveals more similarities between the process $Y_n(t)$ and the BB.

THEOREM 3. *Beside the assumptions of Theorem* 2, *if* $H(x)$ *satisfies*

$$\int \frac{dH(t)}{(1+t\underline{m}(z_1))(1+t\underline{m}(z_2))} - \int \frac{dH(t)}{(1+t\underline{m}(z_1))} \int \frac{dH(t)}{(1+t\underline{m}(z_2))} = 0, \tag{1.13}$$

*then all results of Theorem* 2 *remain true. Moreover, formula* (1.12) *can be simplified to*

$$\begin{aligned}\operatorname{Cov}(X_{g_1}, X_{g_2}) = \frac{2}{c}\biggl(&\int g_1(x)g_2(x)\, dF^{c,H}(x) \\ &- \int g_1(x_1)\, dF^{c,H}(x_1) \int g_2(x_2)\, dF^{c,H}(x_2)\biggr).\end{aligned} \tag{1.14}$$

REMARK 6. Obviously, (1.13) holds when $T_n = bI$. Actually, (1.13) holds if and only if $H(x)$ is a degenerate distribution. To see this, one need only choose $z_2$ to be the complex conjugate of $z_1$.

REMARK 7. Theorem 3 extends the theorem of Silverstein [15]. First, one sees that the $r$th moment of $F_1^{A_n}(x)$ is $x_n^* A_n^r x_n$, which is a special case with $g(x) = x^r$. Then applying Theorem 3 with $T_n = bI$ and combining with Theorem 1.1 of [7], one can obtain the sufficient part of (a) in the theorem of Silverstein [15]. Actually, for $T_n = I$, formula (1.12) can be simplified to

$$\begin{aligned}\operatorname{Cov}(X_{g_1}, X_{g_2}) = \frac{2}{c}\biggl(&\int g_1(x)g_2(x)\, dF_c(x) \\ &- \int g_1(x_1)\, dF_c(x_1) \int g_2(x_2)\, dF_c(x_2)\biggr),\end{aligned} \tag{1.15}$$

where $F_c(x)$ is a special case of $F^{c,H}(x)$, as $T_n = I$.

The organization of the rest of the paper is as follows. In the next section, we complete the proof of Theorem 1. The proof of Theorem 2 is converted to an intermediate Lemma 2, given in Section 3. Sections 4 and 5 contain the proof of Lemma 2. Theorem 3 and some comparisons with the results of [15] are given in Section 6. A truncation lemma (Lemma 4) is postponed to Section 7.

**2. Proof of Theorem 1.** Without loss of generality, we assume that $\|T_n\| \le 1$, where $\|\cdot\|$ denotes the spectral norm on the matrices, that is, their largest singular values. Throughout this paper, $K$ denotes a universal constant which may take different values at different appearances.

LEMMA 1. (Lemma 2.7 in [5]). *Let $X = (X_1, \dots, X_n)$, where $X_j$'s are i.i.d. complex random variables with zero mean and unit variance. Let $B$ be a deterministic $n \times n$ complex matrix. Then for any $p \ge 2$, we have*

$$E|X^*BX - \operatorname{tr} B|^p \le K_p((E|X_1|^4 \operatorname{tr} BB^*)^{p/2} + E|X_1|^{2P} \operatorname{tr}(BB^*)^{p/2}).$$

For $K > 0$, let $\tilde X_{ij} = X_{ij} I(|X_{ij}| \le K) - EX_{ij} I(|X_{ij}| \le K)$ and $\tilde A_n = \frac{1}{N} T_n^{1/2} \tilde X_n \tilde X_n^* T_n^{1/2}$, where $\tilde X_n = (\tilde X_{ij})$. Let $v = \Im z > 0$. Since $X_{ij} - \tilde X_{ij} = X_{ij} I(|X_{ij}| > K) - EX_{ij} I(|X_{ij}| > K)$ and $\|(A_n - zI)^{-1}\|$ is bounded by $\frac{1}{v}$, by Theorem 3.1 in [22], we have

$$\begin{aligned}
&|\mathbf{x}_n^*(A_n - zI)^{-1}\mathbf{x}_n - \mathbf{x}_n^*(\tilde A_n - zI)^{-1}\mathbf{x}_n| \\
&\quad \le \|(A_n - zI)^{-1} - (\tilde A_n - zI)^{-1}\| \\
&\quad \le \|(A_n - zI)^{-1}\| \|(A_n - \bar A_n)\| \|(\tilde A_n - zI)^{-1}\| \\
&\quad \le \frac{1}{Nv^2}(\|X_n - \tilde X_n\| \|X_n^*\| + \|\tilde X_n\| \|X_n^* - \tilde X_n^*\|) \\
&\quad \to \frac{(1+\sqrt{c})^2}{v^2}[E^{1/2}|X_{11} - \tilde X_{11}|^2 (E^{1/2}|\tilde X_{11}|^2 + E^{1/2}|X_{11}|^2)] \qquad \text{a.s.} \\
&\quad \le \frac{2(1+\sqrt{c})^2}{v^2} E^{1/2}|X_{11}|^2 I(|X_{11}| > K).
\end{aligned}$$

The above bound can be made arbitrarily small by choosing $K$ sufficiently large. Since $\lim_{n\to\infty} E|\bar X_{11}|^2 = 1$, the rescaling of $\tilde X_{ij}$ can be dealt with similarly. Hence, in the sequel, it is enough to assume that $|X_{ij}| \le K, EX_{11} = 0$ and $E|X_{11}|^2 = 1$ (for simplicity, suppressing all super- and subscripts on the variables $X_{ij}$).

Next, we will show that

$$\mathbf{x}_n^*(A_n - zI)^{-1}\mathbf{x}_n - \mathbf{x}_n^* E(A_n - zI)^{-1}\mathbf{x}_n \to 0 \qquad \text{a.s.} \tag{2.1}$$

Let $\mathbf{s}_j$ denote the $j$th column of $\frac{1}{\sqrt{N}}T_n^{1/2}X_n$, $A(z)=A_n-zI$, $A_j(z)=A(z)-\mathbf{s}_j\mathbf{s}_j^*$,

$$\alpha_j(z)=\mathbf{s}_j^*A_j^{-1}(z)\mathbf{x}_n\mathbf{x}_n^*(E\underline{m}_n(z)T_n+I)^{-1}\mathbf{s}_j - \frac{1}{N}\mathbf{x}_n^*(E\underline{m}_n(z)T_n+I)^{-1}T_nA_j^{-1}(z)\mathbf{x}_n,$$

$$\xi_j(z)=\mathbf{s}_j^*A_j^{-1}(z)\mathbf{s}_j-\frac{1}{N}\operatorname{tr}T_nA_j^{-1}(z),$$

$$\gamma_j=\mathbf{s}_j^*A_j^{-1}(z)\mathbf{x}_n\mathbf{x}_n^*A_j^{-1}(z)\mathbf{s}_j-\frac{1}{N}\mathbf{x}_n^*A_j^{-1}(z)T_nA_j^{-1}(z)\mathbf{x}_n$$

and

$$\beta_j(z)=\frac{1}{1+\mathbf{s}_j^*A_j^{-1}(z)\mathbf{s}_j},\qquad b_j(z)=\frac{1}{1+N^{-1}\operatorname{tr}T_nA_j^{-1}(z)}.$$

Noting that $|\beta_j(z)|\le|z|/v$ and $\|A_j^{-1}(z)\|\le 1/v$, by Lemma 1, we have

$$E|\mathbf{s}_j^*A_j^{-1}(z)\mathbf{x}_n\mathbf{x}_n^*A_j^{-1}(z)\mathbf{s}_j|^r=O\left(\frac{1}{N^r}\right),\qquad E|\xi_j(z)|^r=O\left(\frac{1}{N^{r/2}}\right). \tag{2.2}$$

Define the $\sigma$-field $\mathcal{F}_j=\sigma(\mathbf{s}_1,\dots,\mathbf{s}_j)$, let $E_j(\cdot)$ denote conditional expectation given the $\sigma$-field $\mathcal{F}_j$ and let $E_0(\cdot)$ denote the unconditional expectation. Note that

$$\begin{aligned}
&\mathbf{x}_n^*(A_n-zI)^{-1}\mathbf{x}_n-\mathbf{x}_n^*E(A_n-zI)^{-1}\mathbf{x}_n\\
&\quad=\sum_{j=1}^N \mathbf{x}_n^*E_jA^{-1}(z)\mathbf{x}_n-\mathbf{x}_n^*E_{j-1}A^{-1}(z)\mathbf{x}_n\\
&\quad=\sum_{j=1}^N \mathbf{x}_n^*E_j(A^{-1}(z)-A_j^{-1}(z))\mathbf{x}_n-\mathbf{x}_n^*E_{j-1}(A^{-1}(z)-A_j^{-1}(z))\mathbf{x}_n\\
&\quad=-\sum_{j=1}^N (E_j-E_{j-1})\beta_j(z)\mathbf{s}_j^*A_j^{-1}(z)\mathbf{x}_n\mathbf{x}_n^*A_j^{-1}(z)\mathbf{s}_j\\
&\quad=-\sum_{j=1}^N E_jb_j(z)\gamma_j(z)\\
&\qquad\qquad-(E_j-E_{j-1})\mathbf{s}_j^*A_j^{-1}(z)\mathbf{x}_n\mathbf{x}_n^*A_j^{-1}(z)\mathbf{s}_j\beta_j(z)b_j(z)\xi_j(z).
\end{aligned}\tag{2.3}$$

By the fact that $|\frac{1}{1+\mathbf{s}_j^*A_j^{-1}(z)\mathbf{s}_j}|\le\frac{|z|}{v}$ and use of the Burkholder inequality, (2.2) and the martingale expression (2.3), we have

$$E|\mathbf{x}_n^*(A_n-zI)^{-1}\mathbf{x}_n-\mathbf{x}_n^*E(A_n-zI)^{-1}\mathbf{x}_n|^r$$

$$\leq E\left[\sum_{j=1}^{N} E_{j-1}|(E_j - E_{j-1})\beta_j(z)\mathbf{s}_j^* A_j^{-1}(z)\mathbf{x}_n\mathbf{x}_n^* A_j^{-1}(z)\mathbf{s}_j|^2\right]^{r/2}$$

$$+ E\sum_{j=1}^{N} |(E_j - E_{j-1})\beta_j(z)\mathbf{s}_j^* A_j^{-1}(z)\mathbf{x}_n\mathbf{x}_n^* A_j^{-1}(z)\mathbf{s}_j|^r$$

$$\leq E\left[\sum_{j=1}^{N} \frac{K|z|^2}{v^2} E_{j-1}|\gamma_j(z)|^2 + E_{j-1}|\mathbf{s}_j^* A_j^{-1}(z)\mathbf{x}_n\mathbf{x}_n^* A_j^{-1}(z)\mathbf{s}_j\xi_j(z)|^2\right]^{r/2}$$

$$+ \sum_{j=1}^{N} \frac{K|z|^r}{v^r} E|\mathbf{s}_j^* A_j^{-1}(z)\mathbf{x}_n\mathbf{x}_n^* A_j^{-1}(z)\mathbf{s}_j|^r$$

$$\leq K[N^{-r/2} + N^{-r+1}].$$

Thus (2.1) follows from Borel–Cantelli lemma, by taking $r > 2$.

Write

$$A(z) - (-zE\underline{m}_n(z)T_n - zI) = \sum_{j=1}^{N} \mathbf{s}_j\mathbf{s}_j^* - (-zE\underline{m}_n(z))T_n.$$

Using the identities

$$\mathbf{s}_j^* A^{-1}(z) = \beta_j(z)\mathbf{s}_j^* A_j^{-1}(z)$$

and

$$\underline{m}_n(z) = -\frac{1}{zN}\sum_{j=1}^{N}\beta_j(z) \tag{2.4}$$

(see (2.2) of [13]), we obtain

$$EA^{-1}(z) - (-zE\underline{m}_n(z)T_n - zI)^{-1}$$

$$= (zE\underline{m}_n(z)T_n + zI)^{-1} E\left[\sum_{j=1}^{N} \mathbf{s}_j\mathbf{s}_j^* - (-zE\underline{m}_n(z))T_n A^{-1}(z)\right]$$

$$= \frac{1}{z}\sum_{j=1}^{N} E\beta_j(z)\Big[(E\underline{m}_n(z)T_n + I)^{-1}\mathbf{s}_j\mathbf{s}_j^* A_j^{-1}(z)$$

$$- \frac{1}{N}(E\underline{m}_n(z)T_n + I)^{-1} T_n EA^{-1}(z)\Big].$$

Multiplying by $\mathbf{x}_n^*$ on the left and $\mathbf{x}_n$ on the right, we have

$$\mathbf{x}_n^* EA^{-1}(z)\mathbf{x}_n - \mathbf{x}_n^*(-zE\underline{m}_n(z)T_n - zI)^{-1}\mathbf{x}_n$$

$$
\begin{aligned}
&= N\frac{1}{z}E\beta_1(z)\bigg[\mathbf{s}_1^* A_1^{-1}(z)\mathbf{x}_n\mathbf{x}_n^*(E\underline{m}_n(z)T_n+I)^{-1}\mathbf{s}_1 \\
&\qquad\qquad - \frac{1}{N}\mathbf{x}_n^*(E\underline{m}_n(z)T_n+I)^{-1}T_nEA^{-1}(z)\mathbf{x}_n\bigg] \qquad (2.5)\\
&\stackrel{\triangle}{=} \delta_1+\delta_2+\delta_3,
\end{aligned}
$$

where

$$
\begin{aligned}
\delta_1 &= \frac{N}{z}E\beta_1(z)\alpha_1(z),\\
\delta_2 &= \frac{1}{z}E\beta_1(z)\mathbf{x}_n^*(E\underline{m}_n(z)T_n+I)^{-1}T_n(A_1^{-1}(z)-A^{-1}(z))\mathbf{x}_n,\\
\delta_3 &= \frac{1}{z}E\beta_1(z)\mathbf{x}_n^*(E\underline{m}_n(z)T_n+I)^{-1}T_n(A^{-1}(z)-EA^{-1}(z))\mathbf{x}_n.
\end{aligned}
$$

Similar to (2.2), by Lemma 1, for $r\ge 2$, we have

$$
E|\alpha_j(z)|^r = O\bigg(\frac{1}{N^r}\bigg).
$$

Therefore,

$$
\delta_1 = -\frac{N}{z}Eb_1(z)\beta_1(z)\xi_1(z)\alpha_1(z) = O(N^{-3/2}). \tag{2.6}
$$

It follows that

$$
\begin{aligned}
|\delta_2| &= \frac{1}{|z|}|E\beta_1^2(z)\mathbf{x}_n^*(E\underline{m}_n(z)T_n+I)^{-1}T_nA_1^{-1}(z)\mathbf{s}_1\mathbf{s}_1^*A_1^{-1}(z)\mathbf{x}_n|\\
&\le K(|E|\mathbf{x}_n^*(E\underline{m}_n(z)T_n+I)^{-1}T_nA_1^{-1}(z)\mathbf{s}_1|^2E|\mathbf{s}_1^*A_1^{-1}(z)\mathbf{x}_n|^2)^{1/2} \qquad (2.7)\\
&= O(N^{-1})
\end{aligned}
$$

and

$$
\begin{aligned}
|\delta_3| &= \frac{1}{|z|}|E\beta_1(z)b_1(z)\xi_1(z)\mathbf{x}_n^*(E\underline{m}_n(z)T_n+I)^{-1}\\
&\qquad\qquad \times T_n(A^{-1}(z)-EA^{-1}(z))\mathbf{x}_n|\\
&\le K(E|\xi_1(z)|^2E|\mathbf{x}_n^*(E\underline{m}_n(z)T_n+I)^{-1} \qquad (2.8)\\
&\qquad\qquad \times T_n(A^{-1}(z)-EA^{-1}(z))\mathbf{x}_n|^2)^{1/2}\\
&= o(N^{-1/2}),
\end{aligned}
$$

where to estimate the second factor, we need to use the martingale decomposition of $A^{-1}(z)-EA^{-1}(z)$.

Combining the above three results and (2.5), we conclude that

$$\mathbf{x}_n^* E A^{-1}(z)\mathbf{x}_n - \mathbf{x}_n^*(-zE\underline{m}_n(z)T_n - zI)^{-1}\mathbf{x}_n \to 0. \tag{2.9}$$

In [13], it is proved that, under the conditions of Theorem 1, $E\underline{m}_n(z) \to \underline{m}(z)$, which is the solution of equation (1.2), and we then conclude that

$$\mathbf{x}_n^* E A^{-1}(z)\mathbf{x}_n - \mathbf{x}_n^*(-z\underline{m}(z)T_n - zI)^{-1}\mathbf{x}_n \to 0.$$

By condition (3) of Theorem 1, we finally obtain that

$$\mathbf{x}_n^* E A^{-1}(z)\mathbf{x}_n \to \int \frac{dH(t)}{-z\underline{m}t - z},$$

which completes the proof of Theorem 1.

**3. An intermediate lemma.** In the sequel, we will follow the work of Bai and Silverstein [7]. To complete the proof of Theorem 2, we need an intermediate lemma.

Write

$$M_n(z) = \sqrt{N}(m_{F_1^{A_n}}(z) - m_{F^{c_n,H_n}}(z)),$$

which is defined on a contour $\mathcal{C}$ in the complex plane, described as follows. Let $u_r$ be a number which is greater than the right endpoint of interval (1.10) and let $u_l$ be a negative number if the left endpoint of interval (1.10) is zero, otherwise let $u_l \in (0, \liminf_n \lambda_{\min}^{T_n} I_{(0,1)}(c)(1-\sqrt{c})^2)$. Let $v_0 > 0$ be arbitrary. Define

$$\mathcal{C}_u = \{u + iv_0 : u \in [u_l, u_r]\}.$$

Then the contour

$$\mathcal{C} = \mathcal{C}_u \cup \{u_l + iv : v \in [0, v_0]\} \cup \{u_r + iv : v \in [0, v_0]\}$$
$$\cup \{\text{their symmetric parts below the real axis}\}.$$

Under the conditions of Theorem 2, for later use, we may select the contour $\mathcal{C}$ in the region on which the functions $g$ are analytic.

As in [7], due to technical difficulties, we will consider $M_n^*(z)$, a truncated version of $M_n(z)$. Choose a sequence of positive numbers $\{\delta_n\}$ such that for $0 < \rho < 1$,

$$\delta_n \downarrow 0, \qquad \delta_n \ge n^{-\rho}. \tag{3.1}$$

Write

$$\mathcal{C}_l = \begin{cases} \{u_l + iv : v \in [n^{-1}\delta_n, v_0]\}, & \text{if } u_l > 0, \\ \{u_l + iv : v \in [0, v_0]\}, & \text{if } u_l < 0, \end{cases}$$

and

$$\mathcal{C}_r = \{u_r + iv : v \in [n^{-1}\delta_n, v_0]\}.$$

Let $\mathcal{C}_0 = \mathcal{C}_u \cup \mathcal{C}_l \cup \mathcal{C}_r$. Now, for $z = u + iv$, we can define the process

$$M_n^*(z) = \begin{cases} M_n(z), \\ \qquad \text{if } z \in \mathcal{C}_0 \cup \bar{\mathcal{C}_0}, \\ \dfrac{nv+\delta_n}{2\delta_n} M_n(u_r + in^{-1}\delta_n) + \dfrac{\delta_n - nv}{2\delta_n} M_n(u_r - in^{-1}\delta_n), \\ \qquad \text{if } u = u_r,\ v \in [-n^{-1}\delta_n, n^{-1}\delta_n], \\ \dfrac{nv+\delta_n}{2\delta_n} M_n(u_l + in^{-1}\delta_n) + \dfrac{\delta_n - nv}{2\delta_n} M_n(u_l - in^{-1}\delta_n), \\ \qquad \text{if } u = u_l > 0,\ v \in [-n^{-1}\delta_n, n^{-1}\delta_n]. \end{cases}$$

$M_n^*(z)$ can be viewed as a random element in the metric space $C(\mathcal{C}, \mathbb{R}^2)$ of continuous functions from $\mathcal{C}$ to $\mathbb{R}^2$. We shall prove the following lemma.

LEMMA 2. *Under the assumptions of Theorem 1 and assumptions* (4) *and* (5) *of Theorem 2, $M_n^*(z)$ forms a tight sequence on $\mathcal{C}$. Furthermore, when the conditions in* (b) *and* (c) *of Theorem 2 on $X_{11}$ hold, for $z \in \mathcal{C}$, $M_n^*(z)$ converges to a Gaussian process $M(\cdot)$ with zero mean and for $z_1, z_2 \in \mathcal{C}$, under the assumptions in* (b),

$$\operatorname{Cov}(M(z_1), M(z_2)) = \frac{2(z_2\underline{m}(z_2) - z_1\underline{m}(z_1))^2}{c^2 z_1 z_2 (z_2 - z_1)(\underline{m}(z_2) - \underline{m}(z_1))}, \tag{3.2}$$

*while under the assumptions in* (c), *the covariance function similar to* (3.2) *is the half of the value of* (3.2).

Similar to the approach of [7], to prove Theorem 2, it suffices to prove Lemma 2. Before proceeding with the detailed proof of the lemma, we need to truncate, recentralize and renormalize the variables $X_{ij}$. However, those procedures are purely technical (and tedious), thus we shall postpone then to the last section of the paper. Now, according to Lemma 4, we further assume that the underlying variables satisfy the following additional conditions:

$$|X_{ij}| \le \varepsilon_n n^{1/4}, \qquad EX_{11} = 0, \qquad E|X_{11}|^2 = 1, \qquad E|X_{11}|^4 < \infty$$

and

$$\begin{cases} E|X_{11}|^4 = 3 + o(1), \\ \qquad \text{if assumption (b) of Theorem 2 is satisfied}, \\ EX_{11}^2 = o(n^{-1/2}),\ E|X_{11}|^4 = 2 + o(1), \\ \qquad \text{if assumption (c) of Theorem 2 is satisfied}. \end{cases}$$

Here, $\varepsilon_n$ is a sequence of positive numbers which converges to zero.

The proof of Lemma 2 will be given in the next two sections.

**4. Convergence in finite dimensions.** For $z\in\mathcal{C}_0$, let

$$M_n^1(z)=\sqrt{N}(m_{F_1^{A_n}}(z)-Em_{F_1^{A_n}}(z))$$

and

$$M_n^2(z)=\sqrt{N}(Em_{F_1^{A_n}}(z)-m_{F^{c_n,H_n}}(z)).$$

Then

$$M_n(z)=M_n^1(z)+M_n^2(z).$$

In this section for any positive integer $r$ and complex numbers $a_1,\dots,a_r$, we will show that

$$\sum_{i=1}^{r}a_iM_n^1(z_i)\qquad(\Im z_i\neq 0)$$

converges in distribution to a Gaussian random variable and will derive the covariance function (3.2). To this end, we employ the notation introduced in Section 2.

Before proceeding with the proofs, we first recall some known facts and results.

1. (See [7].) Let $Y=(Y_1,\dots,Y_n)$, where $Y_i$'s are i.i.d. complex random variables with mean zero and variance 1. Let $A=(a_{ij})_{n\times n}$ and $B=(b_{ij})_{n\times n}$ be complex matrices. Then the following identity holds:

$$\begin{aligned}&E(Y^*AY-\operatorname{tr}A)(Y^*BY-\operatorname{tr}B)\\&\quad=(E|Y_1|^4-|EY_1^2|^2-2)\sum_{i=1}^{n}a_{ii}b_{ii}+|EY_1^2|^2\operatorname{tr}AB^T+\operatorname{tr}AB.\end{aligned}\tag{4.1}$$

2. (See Theorem 35.12 of [9].)

LEMMA 3. *Suppose that for each $n$, $Y_{n1},Y_{n2},\dots,Y_{nr_n}$ is a real martingale difference sequence with respect to the increasing $\sigma$-field $\{\mathcal{F}_{nj}\}$ having second moments. If, as $n\to\infty$, we have*

$$\text{(i)}\qquad \sum_{j=1}^{r_n}E(Y_{nj}^2|\mathcal{F}_{n,j-1})\xrightarrow{i.p.}\sigma^2,$$

*where $\sigma^2$ is a positive constant and for each $\varepsilon>0$,*

$$\text{(ii)}\qquad \sum_{j=1}^{r_n}E(Y_{nj}^2I_{(|Y_{nj}|\geq\varepsilon)})\to 0,$$

*then*

$$\sum_{j=1}^{r_n}Y_{nj}\xrightarrow{\mathcal{D}}N(0,\sigma^2).$$

3. Some simple results follow by using the truncation and centralization steps described in Lemma 4 given in the Appendix:

$$\begin{aligned} E|\mathbf{s}_1^* C \mathbf{s}_1 - N^{-1} \operatorname{tr} T_n C|^p &\le K_p \|C\|^p (\varepsilon_n^{2p-4} N^{-p/2} + N^{-p/2}) \\ &\le K_p \|C\|^p N^{-p/2}, \end{aligned} \tag{4.2}$$

$$E|\mathbf{s}_1^* C \mathbf{x}_n \mathbf{x}_n^* D \mathbf{s}_1 - N^{-1} \mathbf{x}_n^* D T_n C \mathbf{x}_n|^p \le K_p \|C\|^p \|D\|^p \varepsilon_n^{2p-4} N^{-p/2-1}, \tag{4.3}$$

$$E|\mathbf{s}_1^* C \mathbf{x}_n \mathbf{x}_n^* D \mathbf{s}_1|^p \le K_p \|C\|^p \|D\|^p \varepsilon_n^{2p-4} N^{-p/2-1}. \tag{4.4}$$

Let $v = \Im z$. To facilitate the analysis, we will assume that $v > 0$. By (2.3), we have

$$\begin{aligned} &\sqrt{N}(m_{F_1^{A_n}}(z) - E m_{F_1^{A_n}}(z)) \\ &\qquad = -\sqrt{N} \sum_{j=1}^{N} (E_j - E_{j-1}) \beta_j(z) \mathbf{s}_j^* A_j^{-1}(z) \mathbf{x}_n \mathbf{x}_n^* A_j^{-1}(z) \mathbf{s}_j. \end{aligned}$$

Since

$$\beta_j(z) = b_j(z) - \beta_j(z) b_j(z) \xi_j(z) = b_j(z) - b_j^2(z) \xi_j(z) + b_j^2(z) \beta_j(z) \xi_j^2(z),$$

we get

$$\begin{aligned} &(E_j - E_{j-1}) \beta_j(z) \mathbf{s}_j^* A_j^{-1}(z) \mathbf{x}_n \mathbf{x}_n^* A_j^{-1}(z) \mathbf{s}_j \\ &\quad = E_j b_j(z) \gamma_j(z) - E_j \left( b_j^2(z) \xi_j(z) \frac{1}{N} \mathbf{x}_n^* A_j^{-1}(z) T_n A_j^{-1}(z) \mathbf{x}_n \right) \\ &\qquad + (E_j - E_{j-1}) (b_j^2(z) \beta_j(z) \xi_j^2(z) \mathbf{s}_j^* A_j^{-1}(z) \mathbf{x}_n \mathbf{x}_n^* A_j^{-1}(z) \mathbf{s}_j \\ &\qquad\qquad - b_j^2(z) \xi_j(z) \gamma_j(z)). \end{aligned}$$

Applying (4.2), we obtain

$$\begin{aligned} &E \left| \sqrt{N} \sum_{j=1}^{N} E_j \left( b_j^2(z) \xi_j(z) \frac{1}{N} \mathbf{x}_n^* A_j^{-1}(z) T_n A_j^{-1}(z) \mathbf{x}_n \right) \right|^2 \\ &\quad = \frac{1}{N} \sum_{j=1}^{N} E|E_j(b_j^2(z) \xi_j(z) \mathbf{x}_n^* A_j^{-1}(z) T_n A_j^{-1}(z) \mathbf{x}_n)|^2 \\ &\quad \le K \frac{|z|^4}{v^8} E|\xi_1(z)|^2 = O(N^{-1}), \end{aligned}$$

which implies that $\sqrt{N} \sum_{j=1}^{N} E_j(b_j^2(z) \xi_j(z) \frac{1}{N} \mathbf{x}_n^* A_j^{-1}(z) T_n A_j^{-1}(z) \mathbf{x}_n) \xrightarrow{i.p.} 0$.

By (4.2), (4.4) and Hölder's inequality, we have

$$
\begin{aligned}
&E\left|\sqrt{N}\sum_{j=1}^{N}(E_j-E_{j-1})b_j^2(z)\beta_j(z)\xi_j^2(z)\mathbf{s}_j^*A_j^{-1}(z)\mathbf{x}_n\mathbf{x}_n^*A_j^{-1}(z)\mathbf{s}_j\right|^2\\
&\qquad\le K\left(\frac{|z|}{v}\right)^6N^2\left(E|\xi_1^4(z)||\gamma_1^2(z)|+\frac{1}{N^2}E|\xi_1^4(z)||\mathbf{x}_n^*A_1^{-1}(z)T_nA_1^{-1}(z)\mathbf{x}_n|^2\right)\\
&\qquad=O(N^{-3/2}),
\end{aligned}
$$

which implies that

$$
\sqrt{N}\sum_{j=1}^{N}(E_j-E_{j-1})b_j^2(z)\beta_j(z)\xi_j^2(z)\mathbf{s}_j^*A_j^{-1}(z)\mathbf{x}_n\mathbf{x}_n^*A_j^{-1}(z)\mathbf{s}_j\xrightarrow{i.p.}0.
$$

Using a similar argument, we have

$$
\sqrt{N}\sum_{j=1}^{N}(E_j-E_{j-1})b_j^2(z)\xi_j(z)\gamma_j(z)\xrightarrow{i.p.}0.
$$

The estimate above (4.3) of [5], (2.17) of [7] and (4.3) collectively yield that

$$
\begin{aligned}
&E|(b_j(z)+z\underline{m}(z))\gamma_j(z)|^2\\
&\qquad=E[E(|(b_j(z)+z\underline{m}(z))\gamma_j(z)|^2|\sigma(\mathbf{s}_i,i\neq j))]\\
&\qquad=E[|b_j(z)+z\underline{m}(z)|^2E(|\gamma_j(z)|^2|\sigma(\mathbf{s}_i,i\neq j))]=o(N^{-2}),
\end{aligned}
$$

which gives

$$
\sqrt{N}\sum_{j=1}^{N}E_j[(b_j(z)+z\underline{m}(z))\gamma_j(z)]\xrightarrow{i.p.}0.
$$

Note that the above results also hold when $\Im z\le -v_0$, by symmetry. Hence for finite-dimensional convergence, we need only consider the sum

$$
\sum_{i=1}^{r}a_i\sum_{j=1}^{N}Y_j(z_i)=\sum_{j=1}^{N}\sum_{i=1}^{r}a_iY_j(z_i),
$$

where $Y_j(z_i)=-\sqrt{N}z_i\underline{m}(z_i)E_j\gamma_j(z_i)$. Since

$$
E|Y_j(z_i)|^4=O(\varepsilon_n^4N^{-1}),
$$

we have

$$
\sum_{j=1}^{N}E\left(\left|\sum_{i=1}^{r}a_iY_j(z_i)\right|^2I\left(\left|\sum_{i=1}^{r}a_iY_j(z_i)\right|\ge\varepsilon\right)\right)\le\frac{1}{\varepsilon^2}\sum_{j=1}^{N}E\left|\sum_{i=1}^{r}a_iY_j(z_i)\right|^4=O(\varepsilon_n^4).
$$

Thus condition (ii) of Lemma 3 is satisfied.

Now, we only need to show that for $z_1, z_2 \in \mathbb{C} \setminus \mathbb{R}$,

$$\sum_{j=1}^{N} E_{j-1}(Y_j(z_1)Y_j(z_2)) \tag{4.5}$$

converges in probability to a constant under the assumptions in (b) or (c). It is easy to verify that

$$|\operatorname{tr} E_j(A_j^{-1}(z_1)\mathbf{x}_n\mathbf{x}_n^* A_j^{-1}(z_1))T_n E_j(A_j^{-1}(z_2)\mathbf{x}_n\mathbf{x}_n^* A_j^{-1}(z_2)T_n)| \le \frac{1}{|v_1 v_2|^2}, \tag{4.6}$$

where $v_1 = \Im(z_1)$ and $v_2 = \Im(z_2)$. It follows that, for the complex case, applying (4.1), (4.5) now becomes

$$\begin{aligned} &z_1 z_2 \underline{m}(z_1)\underline{m}(z_2)\frac{1}{N}\sum_{j=1}^{N} E_{j-1} \operatorname{tr} E_j(A_j^{-1}(z_1)\mathbf{x}_n\mathbf{x}_n^* A_j^{-1}(z_1))T_n \\ &\qquad \times E_j(A_j^{-1}(z_2)\mathbf{x}_n\mathbf{x}_n^* A_j^{-1}(z_2)T_n) + o_p(1) \\ &= z_1 z_2 \underline{m}(z_1)\underline{m}(z_2)\frac{1}{N}\sum_{j=1}^{N} E_{j-1}(\mathbf{x}_n^* A_j^{-1}(z_1)T_n \breve{A}_j^{-1}(z_2)\mathbf{x}_n) \\ &\qquad \times (\mathbf{x}_n^* \breve{A}_j^{-1}(z_2)T_n A_j^{-1}(z_1)\mathbf{x}_n) + o_p(1), \end{aligned} \tag{4.7}$$

where $\breve{A}_j^{-1}(z_2)$ is defined similarly as $A_j^{-1}(z_2)$ by $(\mathbf{s}_1, \ldots, \mathbf{s}_{j-1}, \breve{\mathbf{s}}_{j+1}, \ldots, \breve{\mathbf{s}}_N)$ and where $\breve{\mathbf{s}}_{j+1}, \ldots, \breve{\mathbf{s}}_N$ are i.i.d. copies of $\mathbf{s}_{j+1}, \ldots, \mathbf{s}_N$.

For the real case, (4.5) will be twice the magnitude of (4.7).

Write

$$\begin{aligned} &\mathbf{x}_n^*(A_j^{-1}(z_1) - E_{j-1}A_j^{-1}(z_1))T_n \breve{A}_j^{-1}(z_2)\mathbf{x}_n \\ &\qquad = \sum_{t=j}^{N} \mathbf{x}_n^*(E_t A_j^{-1}(z_1) - E_{t-1}A_j^{-1}(z_1))T_n \breve{A}_j^{-1}(z_2)\mathbf{x}_n. \end{aligned} \tag{4.8}$$

From (4.8), we note that

$$\begin{aligned} &E_{j-1}(\mathbf{x}_n^* A_j^{-1}(z_1)T_n \breve{A}_j^{-1}(z_2)\mathbf{x}_n)(\mathbf{x}_n^* \breve{A}_j^{-1}(z_2)T_n A_j^{-1}(z_1)\mathbf{x}_n) \\ &\quad = \sum_{t=j}^{N} E_{j-1}\mathbf{x}_n^*(E_t A_j^{-1}(z_1) - E_{t-1}A_j^{-1}(z_1))T_n \breve{A}_j^{-1}(z_2)\mathbf{x}_n\mathbf{x}_n^* \breve{A}_j^{-1}(z_2) \\ &\qquad \times T_n(E_t A_j^{-1}(z_1) - E_{t-1}A_j^{-1}(z_1))\mathbf{x}_n \\ &\qquad + E_{j-1}(\mathbf{x}_n^*(E_{j-1}A_j^{-1}(z_1)T_n)\breve{A}_j^{-1}(z_2)\mathbf{x}_n) \\ &\qquad \times (\mathbf{x}_n^* \breve{A}_j^{-1}(z_2)T_n(E_{j-1}A_j^{-1}(z_1))\mathbf{x}_n) \end{aligned} \tag{4.9}$$

$$= E_{j-1}(\mathbf{x}_n^*(E_{j-1}A_j^{-1}(z_1)T_n)\breve{A}_j^{-1}(z_2)\mathbf{x}_n)$$
$$\times (\mathbf{x}_n^*\breve{A}_j^{-1}(z_2)T_n(E_{j-1}A_j^{-1}(z_1))\mathbf{x}_n) + O(N^{-1}),$$

where we have used the fact that

$$|E_{j-1}\mathbf{x}_n^*(E_tA_j^{-1}(z_1) - E_{t-1}A_j^{-1}(z_1))$$
$$\times T_n\breve{A}_j^{-1}(z_2)\mathbf{x}_n\mathbf{x}_n^*\breve{A}_j^{-1}(z_2)T_n(E_tA_j^{-1}(z_1) - E_{t-1}A_j^{-1}(z_1))\mathbf{x}_n|$$
$$\le 4(|E_{j-1}|\beta_{tj}(z_1)\mathbf{x}_n^*(A_{tj}^{-1}(z_1)\mathbf{s}_t\mathbf{s}_t^*(A_{tj}^{-1}(z_1)T_n\breve{A}_j^{-1}(z_2)\mathbf{x}_n|^2$$
$$\times E_{j-1}|\beta_{tj}(z_1)\mathbf{x}_n^*\breve{A}_j^{-1}(z_2)T_n(A_{tj}^{-1}(z_1)\mathbf{s}_t\mathbf{s}_t^*A_j^{-1}(z_1))\mathbf{x}_n|^2)^{1/2} = O(N^{-2}).$$

Similarly, one can prove that

$$E_{j-1}(\mathbf{x}_n^*(E_{j-1}A_j^{-1}(z_1)T_n)\breve{A}_j^{-1}(z_2)\mathbf{x}_n)(\mathbf{x}_n^*\breve{A}_j^{-1}(z_2)T_n(E_{j-1}A_j^{-1}(z_1))\mathbf{x}_n)$$
$$= E_{j-1}(\mathbf{x}_n^*A_j^{-1}(z_1)T_n\breve{A}_j^{-1}(z_2)\mathbf{x}_n)$$
$$\times E_{j-1}(\mathbf{x}_n^*\breve{A}_j^{-1}(z_2)T_nA_j^{-1}(z_1)\mathbf{x}_n) + O(N^{-1}).$$

Define

$$A_{ij}(z) = A(z) - \mathbf{s}_i\mathbf{s}_i^* - \mathbf{s}_j\mathbf{s}_j^*, \qquad T^{-1}(z_1) = \left(z_1I - \frac{N-1}{N}b_{n1}(z_1)T_n\right)^{-1},$$
$$\beta_{ij}(z) = \frac{1}{1+\mathbf{s}_i^*A_{ij}^{-1}(z)\mathbf{s}_i} \quad \text{and} \quad b_{n1}(z) = \frac{1}{1+N^{-1}E\operatorname{tr}T_nA_{12}^{-1}(z)}.$$

Then (see (2.9) in [7])

$$A_j(z_1) = -T^{-1}(z_1) + b_{n1}(z_1)B_j(z_1) + C_j(z_1) + D_j(z_1), \tag{4.10}$$

where $B_1(z_1) = B_{j1}(z_1) + B_{j2}(z_1)$,

$$B_{j1}(z_1) = \sum_{i>j} T^{-1}(z_1)(\mathbf{s}_i\mathbf{s}_i^* - N^{-1}T_n)A_{ij}^{-1}(z_1),$$
$$B_{j2}(z_1) = \sum_{i<j} T^{-1}(z_1)(\mathbf{s}_i\mathbf{s}_i^* - N^{-1}T_n)A_{ij}^{-1}(z_1),$$
$$C_j(z_1) = \sum_{i\ne j}(\beta_{ij}(z_1) - b_{n1}(z_1))T^{-1}(z_1)\mathbf{s}_i\mathbf{s}_i^*A_{ij}^{-1}(z_1)$$

and

$$D_j(z_1) = N^{-1}b_{n1}(z_1)T^{-1}(z_1)T_n\sum_{i\ne j}(A_{ij}^{-1}(z_1) - A_j(z_1)).$$

It follows that

$$
\begin{aligned}
&E_{j-1}(\mathbf{x}_n^* A_j^{-1}(z_1) T_n \breve{A}_j^{-1}(z_2)\mathbf{x}_n) E_{j-1}(\mathbf{x}_n^* \breve{A}_j^{-1}(z_2) T_n A_j^{-1}(z_1)\mathbf{x}_n)\\
&\quad = -E_{j-1}(\mathbf{x}_n^* A_j^{-1}(z_1) T_n \breve{A}_j^{-1}(z_2)\mathbf{x}_n)\\
&\qquad \times E_{j-1}(\mathbf{x}_n^* \breve{A}_j^{-1}(z_2) T_n T^{-1}(z_1) T_n \mathbf{x}_n)\\
&\qquad + B(z_1,z_2) + C(z_1,z_2) + D(z_1,z_2),
\end{aligned} \tag{4.11}
$$

where

$$
\begin{aligned}
B(z_1,z_2) &= b_{n1}(z_1) E_{j-1}(\mathbf{x}_n^* A_j^{-1}(z_1) T_n \breve{A}_j^{-1}(z_2)\mathbf{x}_n)\\
&\qquad \times E_{j-1}(\mathbf{x}_n^* \breve{A}_j^{-1}(z_2) T_n B_{j2}(z_1)\mathbf{x}_n),\\
C(z_1,z_2) &= E_{j-1}(\mathbf{x}_n^* A_j^{-1}(z_1) T_n \breve{A}_j^{-1}(z_2)\mathbf{x}_n)\\
&\qquad \times E_{j-1}(\mathbf{x}_n^* \breve{A}_j^{-1}(z_2) T_n C_j(z_1)\mathbf{x}_n)
\end{aligned}
$$

and

$$
D(z_1,z_2) = E_{j-1}(\mathbf{x}_n^* A_j^{-1}(z_1) T_n \breve{A}_j^{-1}(z_2)\mathbf{x}_n)(\mathbf{x}_n^* \breve{A}_j^{-1}(z_2) T_n D_j(z_1)\mathbf{x}_n).
$$

We then prove that

$$
E|C(z_1,z_2)| = o(1) \quad \text{and} \quad E|D(z_1,z_2)| = o(1). \tag{4.12}
$$

Note that although $C$ and $D$ depend on $j$ implicitly, $E|C(z_1,z_2)|$ and $E|D(z_1,z_2)|$ are independent of $j$ since the entries of $X_n$ are i.i.d.

We then have

$$
\begin{aligned}
E|C(z_1,z_2)| &\le \frac{1}{|v_1 v_2|} E|\mathbf{x}_n^* \breve{A}_j^{-1}(z_2) T_n C_j(z_1)\mathbf{x}_n|\\
&\le \frac{1}{|v_1 v_2|} \sum_{i\ne j} (E|\beta_{ij}(z_1) - b_{n1}(z_1)|^2\\
&\qquad \times E|\mathbf{s}_i^* A_{ij}^{-1}(z_1)\mathbf{x}_n \mathbf{x}_n^* (\breve{A}_j^{-1}(z_2)) T_n T^{-1}(z_1)\mathbf{s}_i|^2)^{1/2}.
\end{aligned}
$$

When $i > j$, $\mathbf{s}_i$ is independent of $\breve{A}_j^{-1}(z_2)$. As in the proof of (2.2), we have

$$
E|\mathbf{s}_i^* A_{ij}^{-1}(z_1)\mathbf{x}_n \mathbf{x}_n^* (\breve{A}_j^{-1}(z_2)) T_n T^{-1}(z_1)\mathbf{s}_i|^2 = O(N^{-2}). \tag{4.13}
$$

When $i < j$, by substituting $\breve{A}_j^{-1}(z_2)$ for $\breve{A}_{ij}^{-1}(z_2) - \breve{\beta}_{ij}(z_2)\breve{A}_{ij}^{-1}(z_2)\mathbf{s}_i\mathbf{s}_i^* \times \breve{A}_{ij}^{-1}(z_2)$, we can also obtain the above inequality. Noting that

$$
E|\beta_{ij}(z_1) - b_{n1}(z_1)|^2 = E|\beta_{ij}(z_1) b_{n1}(z_1)\xi_{ij}|^2 = O(n^{-1}), \tag{4.14}
$$

where $\xi_{ij}(z) = \mathbf{s}_i^* A_{ij}^{-1}(z)\mathbf{s}_i - \frac{1}{N}\operatorname{tr} A_{ij}^{-1}(z)$ and $\breve{\beta}_{ij}(z_2)$ is defined similarly to $\beta_{ij}(z_2)$, and combining (4.13)–(4.14), we conclude that

$$
E|C(z_1,z_2)| = o(1).
$$

The argument for $D(z_1,z_2)$ is similar to that of $C(z_1,z_2)$, only simpler, and is therefore omitted. Hence, (4.12) holds.

Next, write

$$B(z_1,z_2)=B_1(z_1,z_2)+B_2(z_1,z_2)+B_3(z_1,z_2), \tag{4.15}$$

where

$$\begin{aligned}B_1(z_1,z_2)&=\sum_{i<j}b_{n1}(z_1)E_{j-1}\mathbf{x}_n^*\beta_{ij}(z_1)A_{ij}^{-1}(z_1)\mathbf{s}_i\mathbf{s}_i^*A_{ij}^{-1}(z_1)T_n\breve{A}_j^{-1}(z_2)\mathbf{x}_n\\&\qquad\times E_{j-1}\mathbf{x}_n^*\breve{A}_j^{-1}(z_2)T_nT^{-1}(z_1)(\mathbf{s}_i\mathbf{s}_i^*-N^{-1}T_n)A_{ij}^{-1}(z_1)\mathbf{x}_n,\\B_2(z_1,z_2)&=\sum_{i<j}b_{n1}(z_1)E_{j-1}\mathbf{x}_n^*A_{ij}^{-1}(z_1)T_n\breve{A}_{ij}^{-1}(z_2)\mathbf{s}_i\mathbf{s}_i^*\breve{A}_{ij}^{-1}(z_2)\breve{\beta}_{ij}(z_2)\mathbf{x}_n\\&\qquad\times E_{j-1}\mathbf{x}_n^*\breve{A}_j^{-1}(z_2)T_nT^{-1}(z_1)(\mathbf{s}_i\mathbf{s}_i^*-N^{-1}T_n)A_{ij}^{-1}(z_1)\mathbf{x}_n\end{aligned}$$

and

$$\begin{aligned}B_3(z_1,z_2)&=\sum_{i<j}b_{n1}(z_1)E_{j-1}\mathbf{x}_n^*A_{ij}^{-1}(z_1)T_n\breve{A}_{ij}^{-1}(z_2)\mathbf{x}_n\\&\qquad\times E_{j-1}\mathbf{x}_n^*\breve{A}_j^{-1}(z_2)T_nT^{-1}(z_1)(\mathbf{s}_i\mathbf{s}_i^*-N^{-1}T_n)A_{ij}^{-1}(z_1)\mathbf{x}_n.\end{aligned}$$

Splitting $\breve{A}_j^{-1}(z_2)$ into the sum of $\breve{A}_{ij}^{-1}(z_2)$ and $-\breve{\beta}_{ij}(z_2)\breve{A}_{ij}^{-1}(z_2)\mathbf{s}_i\mathbf{s}_i^*\breve{A}_{ij}^{-1}(z_2)$, as in the proof of (4.14), one can show that

$$\begin{aligned}&E|B_1(z_1,z_2)|\\&\quad\le\sum_{i<j}|b_{n1}(z_1)|(E|\mathbf{x}_n^*\beta_{ij}(z_1)A_{ij}^{-1}(z_1)\mathbf{s}_i\mathbf{s}_i^*A_{ij}^{-1}(z_1)T_n\breve{A}_j^{-1}(z_2)\mathbf{x}_n|^2\\&\qquad\qquad\times E|\mathbf{x}_n^*\breve{A}_j^{-1}(z_2)T_nT^{-1}(z_1)(\mathbf{s}_i\mathbf{s}_i^*-N^{-1}T_n)A_{ij}^{-1}(z_1)\mathbf{x}_n|^2)^{1/2}\\&\quad=O(N^{-1/2}).\end{aligned}$$

By the same argument, we have

$$E|B_2(z_1,z_2)|=O(N^{-1}).$$

To deal with $B_3(z_1,z_2)$, we again split $\breve{A}_j^{-1}(z_2)$ into the sum of $\breve{A}_{ij}^{-1}(z_2)$ and $-\breve{\beta}_{ij}(z_2)\breve{A}_{ij}^{-1}(z_2)\mathbf{s}_i\mathbf{s}_i^*\breve{A}_{ij}^{-1}(z_2)$. We first show that

$$\begin{aligned}B_{31}(z_1,z_2)&=\sum_{i<j}b_{n1}(z_1)E_{j-1}\mathbf{x}_n^*A_{ij}^{-1}(z_1)T_n\breve{A}_{ij}^{-1}(z_2)\mathbf{x}_n\\&\qquad\times E_{j-1}\mathbf{x}_n^*\breve{A}_{ij}^{-1}(z_2)T_nT^{-1}(z_1)(\mathbf{s}_i\mathbf{s}_i^*-N^{-1}T_n)A_{ij}^{-1}(z_1)\mathbf{x}_n\\&=o_p(1).\end{aligned}\tag{4.16}$$

We have

$$E|B_{31}(z_1,z_2)|^2=\sum_{i_1,i_2<j}|b_{n1}(z_1)|^2EE_{j-1}\mathbf{x}_n^*A_{i_1j}^{-1}(z_1)T_n\breve{A}_{i_1j}^{-1}(z_2)\mathbf{x}_n$$
$$\times E_{j-1}\mathbf{x}_n^*A_{i_2j}^{-1}(\bar z_1)T_n\breve{A}_{i_2j}^{-1}(\bar z_2)\mathbf{x}_n$$
$$\times E_{j-1}\mathbf{x}_n^*\breve{A}_{i_1j}^{-1}(z_2)T_nT^{-1}(z_1)$$
$$\times(\mathbf{s}_{i_1}\mathbf{s}_{i_1}^*-N^{-1}T_n)A_{i_1j}^{-1}(z_1)\mathbf{x}_n\mathbf{x}_n^*\breve{A}_{i_2j}^{-1}(\bar z_2)$$
$$\times T_nT^{-1}(\bar z_1)(\mathbf{s}_{i_2}\mathbf{s}_{i_2}^*-N^{-1}T_n)A_{i_2j}^{-1}(\bar z_1)\mathbf{x}_n.$$

When $i_1=i_2$, the term in the above expression is bounded by

$$KE|\mathbf{x}_n^*\breve{A}_{i_1j}^{-1}(z_2)T_nT^{-1}(z_1)(\mathbf{s}_{i_1}\mathbf{s}_{i_1}^*-N^{-1}T_n)A_{ij}^{-1}(z_1)\mathbf{x}_n|^2=O(N^{-2}).$$

For $i_1\neq i_2<j$, define

$$\beta_{i_1i_2j}(z_1)=\frac{1}{1+\mathbf{s}_{i_2}^*A_{i_1i_2j}^{-1}(z_1)\mathbf{s}_{i_2}},$$
$$A_{i_1i_2j}(z_1)=A(z_1)-\mathbf{s}_{i_1}\mathbf{s}_{i_1}^*-\mathbf{s}_{i_2}\mathbf{s}_{i_2}^*-\mathbf{s}_j\mathbf{s}_j^*$$

and similarly define $\breve\beta_{i_1,i_2,j}(z_2)$ and $\breve{A}_{i_1i_2j}(z_2)$.

We have

$$|EE_{j-1}\mathbf{x}_n^*A_{i_1,i_2,j}^{-1}(z_1)\mathbf{s}_{i_2}\mathbf{s}_{i_2}^*A_{i_1,i_2,j}^{-1}(z_1)\beta_{i_1,i_2,j}(z_1)T_n\breve{A}_{i_1j}^{-1}(z_2)$$
$$\times\mathbf{x}_nE_{j-1}\mathbf{x}_n^*A_{i_2j}^{-1}(\bar z_1)T_n\breve{A}_{i_2j}^{-1}(\bar z_2)\mathbf{x}_n$$
$$\times E_{j-1}\mathbf{x}_n^*\breve{A}_{i_1j}^{-1}(z_2)T_nT^{-1}(z_1)(\mathbf{s}_{i_1}\mathbf{s}_{i_1}^*-N^{-1}T_n)A_{i_1j}^{-1}(z_1)\mathbf{x}_n$$
$$\times\mathbf{x}_n^*\breve{A}_{i_2j}^{-1}(\bar z_2)T_nT^{-1}(\bar z_1)(\mathbf{s}_{i_2}\mathbf{s}_{i_2}^*-N^{-1}T_n)A_{i_2j}^{-1}(\bar z_1)\mathbf{x}_n|$$
$$\leq K(E|\mathbf{x}_n^*A_{i_1,i_2,j}^{-1}(z_1)\mathbf{s}_{i_2}\mathbf{s}_{i_2}^*A_{i_1,i_2,j}^{-1}(z_1)\beta_{i_1,i_2,j}(z_1)T_n\breve{A}_{i_1j}^{-1}(z_2)\mathbf{x}_n|^2)^{1/2}$$
$$\times(E|\mathbf{x}_n^*\breve{A}_{i_1j}^{-1}(z_2)T_nT^{-1}(z_1)(\mathbf{s}_{i_1}\mathbf{s}_{i_1}^*-N^{-1}T_n)A_{i_1j}^{-1}(z_1)\mathbf{x}_n|^4)^{1/4}$$
$$\times(E|\mathbf{x}_n^*\breve{A}_{i_2j}^{-1}(z_2)T_nT^{-1}(z_1)(\mathbf{s}_{i_2}\mathbf{s}_{i_2}^*-N^{-1}T_n)A_{i_2j}^{-1}(z_1)\mathbf{x}_n|^4)^{1/4}$$
$$=O(N^{-5/2}),$$

$$|EE_{j-1}\mathbf{x}_n^*A_{i_1,i_2,j}^{-1}(z_1)T_n\breve{A}_{i_1,i_2,j}^{-1}(z_2)\mathbf{s}_{i_2}\mathbf{s}_{i_2}^*\breve{A}_{i_1,i_2,j}^{-1}(z_2)\breve\beta_{i_1,i_2,j}(z_2)\mathbf{x}_n$$
$$\times E_{j-1}\mathbf{x}_n^*A_{i_2j}^{-1}(\bar z_1)T_n\breve{A}_{i_2j}^{-1}(\bar z_2)\mathbf{x}_n$$
$$\times E_{j-1}\mathbf{x}_n^*\breve{A}_{i_1j}^{-1}(z_2)T_nT^{-1}(z_1)(\mathbf{s}_{i_1}\mathbf{s}_{i_1}^*-N^{-1}T_n)A_{i_1j}^{-1}(z_1)\mathbf{x}_n$$
$$\times\mathbf{x}_n^*\breve{A}_{i_2j}^{-1}(\bar z_2)T_nT^{-1}(\bar z_1)(\mathbf{s}_{i_2}\mathbf{s}_{i_2}^*-N^{-1}T_n)A_{i_2j}^{-1}(\bar z_1)\mathbf{x}_n|$$

$$\leq K(E|\mathbf{x}_n^* A_{i_1,i_2,j}^{-1}(z_1)T_n\breve{A}_{i_1,i_2,j}^{-1}(z_2)\mathbf{s}_{i_2}\mathbf{s}_{i_2}^*\breve{A}_{i_1,i_2,j}^{-1}(z_2)\breve{\beta}_{i_1,i_2,j}(z_2)\mathbf{x}_n|^2)^{1/2}$$
$$\times (E|\mathbf{x}_n^*\breve{A}_{i_1j}^{-1}(z_2)T_nT^{-1}(z_1)(\mathbf{s}_{i_1}\mathbf{s}_{i_1}^* - N^{-1}T_n)A_{i_1j}^{-1}(z_1)\mathbf{x}_n|^4)^{1/4}$$
$$\times (E|\mathbf{x}_n^*\breve{A}_{i_2j}^{-1}(z_2)T_nT^{-1}(z_1)(\mathbf{s}_{i_2}\mathbf{s}_{i_2}^* - N^{-1}T_n)A_{i_2j}^{-1}(z_1)\mathbf{x}_n|^4)^{1/4}$$
$$= O(N^{-5/2})$$

and, by (4.1),

$$|EE_{j-1}\mathbf{x}_n^* A_{i_1,i_2,j}^{-1}(z_1)T_n\breve{A}_{i_1,i_2,j}^{-1}(z_2)\mathbf{x}_n E_{j-1}\mathbf{x}_n^* A_{i_2j}^{-1}(\bar{z}_1)T_n\breve{A}_{i_2j}^{-1}(\bar{z}_2)\mathbf{x}_n$$
$$\times E_{j-1}\mathbf{x}_n^*\breve{A}_{i_1,i_2,j}^{-1}(z_2)\mathbf{s}_{i_2}\mathbf{s}_{i_2}^*\breve{A}_{i_1,i_2,j}^{-1}(z_2)\breve{\beta}_{i_1,i_2,j}(z_2)T_nT^{-1}(z_1)$$
$$\times (\mathbf{s}_{i_1}\mathbf{s}_{i_1}^* - N^{-1}T_n)$$
$$\times A_{i_1j}^{-1}(z_1)\mathbf{x}_n\mathbf{x}_n^*\breve{A}_{i_2j}^{-1}(\bar{z}_2)T_nT^{-1}(\bar{z}_1)(\mathbf{s}_{i_2}\mathbf{s}_{i_2}^* - N^{-1}T_n)A_{i_2j}^{-1}(\bar{z}_1)\mathbf{x}_n|$$
$$\leq K(E|\mathbf{x}_n^*\breve{A}_{i_1,i_2,j}^{-1}(z_2)\mathbf{s}_{i_2}\mathbf{s}_{i_2}^*\breve{A}_{i_1,i_2,j}^{-1}(z_2)\breve{\beta}_{i_1,i_2,j}(z_2)T_nT^{-1}(z_1)$$
$$\times (\mathbf{s}_{i_1}\mathbf{s}_{i_1}^* - N^{-1}T_n)A_{i_1j}^{-1}(z_1)\mathbf{x}_n|^2)^{1/2}$$
$$\times (E|\mathbf{x}_n^*\breve{A}_{i_2j}^{-1}(z_2)T_nT^{-1}(z_1)(\mathbf{s}_{i_2}\mathbf{s}_{i_2}^* - N^{-1}T_n)A_{i_2j}^{-1}(z_1)\mathbf{x}_n|^2)^{1/2}$$
$$\leq K(E|\mathbf{x}_n^*\breve{A}_{i_1,i_2,j}^{-1}(z_2)\mathbf{s}_{i_2}\mathbf{s}_{i_2}^*\breve{A}_{i_1,i_2,j}^{-1}(z_2)T_nT^{-1}(z_1)$$
$$\times (\mathbf{s}_{i_1}\mathbf{s}_{i_1}^* - N^{-1}T_n)A_{i_1j}^{-1}(z_1)\mathbf{x}_n|^2)^{1/2} \times O(N^{-1})$$
$$\leq K(E\mathbf{x}_n^*\breve{A}_{i_1,i_2,j}^{-1}(z_2)\mathbf{s}_{i_2}\mathbf{s}_{i_2}^*\breve{A}_{i_1,i_2,j}^{-1}(z_2)T_nT^{-1}(z_1)T_nT^{-1}(\bar{z}_1)T_n$$
$$\times \breve{A}_{i_1,i_2,j}^{-1}(\bar{z}_2)\mathbf{s}_{i_2}\mathbf{s}_{i_2}^*\breve{A}_{i_1,i_2,j}^{-1}(\bar{z}_2)\mathbf{x}_n\mathbf{x}_n^* A_{i_1,j}^{-1}(\bar{z}_1)T_nA_{i_1j}^{-1}(z_1)\mathbf{x}_n)^{1/2}$$
$$\times O(N^{-2})$$
$$\leq K(E|\mathbf{x}_n^*\breve{A}_{i_1,i_2,j}^{-1}(z_2)\mathbf{s}_{i_2}\mathbf{s}_{i_2}^*\breve{A}_{i_1,i_2,j}^{-1}(\bar{z}_2)\mathbf{x}_n|^2)^{1/4}$$
$$\times (E|\mathbf{s}_{i_2}^*\breve{A}_{i_1,i_2,j}^{-1}(z_2)T_nT^{-1}(z_1)T_nT^{-1}(\bar{z}_1)T_n\breve{A}_{i_1,i_2,j}^{-1}(\bar{z}_2)\mathbf{s}_{i_2}^2)^{1/4}$$
$$\times O(N^{-2})$$
$$= O(N^{-9/4}).$$

The conclusion (4.16) then follows from the above three estimates.

Therefore,

$$B_3(z_1,z_2) = B_{32}(z_1,z_2) + o_p(1),$$
$$B_{32}(z_1,z_2) = -\sum_{i<j} b_{n1}(z_1)E_{j-1}\mathbf{x}_n^* A_{ij}^{-1}(z_1)T_n\breve{A}_{ij}^{-1}(z_2)\mathbf{x}_n$$
$$\times E_{j-1}\mathbf{x}_n^*\breve{A}_{ij}^{-1}(z_2)\mathbf{s}_i\mathbf{s}_i^*\breve{A}_{ij}^{-1}(z_2)\breve{\beta}_{ij}(z_2)T_nT^{-1}(z_1)$$

$$
\begin{aligned}
&\times (\mathbf{s}_i\mathbf{s}_i^* - N^{-1}T_n)A_{ij}^{-1}(z_1)\mathbf{x}_n \\
&= -\sum_{i<j} b_{n1}(z_1)E_{j-1}\mathbf{x}_n^* A_{ij}^{-1}(z_1)T_n\breve{A}_{ij}^{-1}(z_2)\mathbf{x}_n \\
&\qquad \times E_{j-1}\mathbf{x}_n^*\breve{A}_{ij}^{-1}(z_2)\mathbf{s}_i\mathbf{s}_i^*\breve{A}_{ij}^{-1}(z_2)\breve{\beta}_{ij}(z_2)T_nT^{-1}(z_1) \\
&\qquad \times \mathbf{s}_i\mathbf{s}_i^*A_{ij}^{-1}(z_1)\mathbf{x}_n + o_p(1).
\end{aligned}
$$

By (4.3) of [5], (4.2) and (4.4), for $i<j$, we have

$$
\begin{aligned}
&E|\mathbf{x}_n^*\breve{A}_{ij}^{-1}(z_2)\mathbf{s}_i\mathbf{s}_i^*A_{ij}^{-1}(z_1)\mathbf{x}_n(\mathbf{s}_i^*\breve{A}_{ij}^{-1}(z_2)\breve{\beta}_{ij}(z_2)T_nT^{-1}(z_1)\mathbf{s}_i \\
&\qquad\qquad - N^{-1}b_{n1}(z_2)\operatorname{tr}T_n\breve{A}_{ij}^{-1}(z_2)T_nT^{-1}(z_1))| \\
&\quad\le (E|\mathbf{x}_n^*\breve{A}_{ij}^{-1}(z_2)\mathbf{s}_i\mathbf{s}_i^*A_{ij}^{-1}(z_1)\mathbf{x}_n^2)^{1/2} \\
&\qquad \times [(E|\breve{\beta}_{ij}(z_2)|^2|\mathbf{s}_i^*\breve{A}_{ij}^{-1}(z_2)T_nT^{-1}(z_1)\mathbf{s}_i \\
&\qquad\qquad - N^{-1}\operatorname{tr}T_n\breve{A}_{ij}^{-1}(z_2)T_nT^{-1}(z_1)|^2)^{1/2} \\
&\qquad\quad + (E|\breve{\beta}_{ij}(z_2)-b_{n1}(z_2)|^2|N^{-1}\operatorname{tr}T_n\breve{A}_{ij}^{-1}(z_2)T_nT^{-1}(z_1)|^2)^{1/2}] \\
&\quad = O(N^{-3/2}).
\end{aligned}
\tag{4.17}
$$

Collecting the proofs from (4.10) to (4.17), we have proved that

$$
\begin{aligned}
B(z_1,z_2) &= -b_{n1}(z_1)b_{n1}(z_2) \\
&\quad \times \sum_{i<j} E_{j-1}\mathbf{x}_n^*A_{ij}^{-1}(z_1)T_n\breve{A}_{ij}^{-1}(z_2)\mathbf{x}_nE_{j-1} \\
&\qquad \times (\mathbf{x}_n^*\breve{A}_{ij}^{-1}(z_2)\mathbf{s}_i\mathbf{s}_i^*A_{ij}^{-1}(z_1)\mathbf{x}_nN^{-1}\operatorname{tr}T_n\operatorname{tr}T_n\breve{A}_{ij}^{-1}(z_2)T_nT^{-1}(z_1)) \\
&\quad + o_p(1).
\end{aligned}
$$

Similarly to the proof of (4.16), we may further replace $\mathbf{s}_i\mathbf{s}_i^*$ in the above expression by $N^{-1}T_n$, that is,

$$
\begin{aligned}
B(z_1,z_2) &= -b_{n1}(z_1)b_{n1}(z_2)N^{-2} \\
&\quad \times \sum_{i<j} E_{j-1}\mathbf{x}_n^*A_{ij}^{-1}(z_1)T_n\breve{A}_{ij}^{-1}(z_2)\mathbf{x}_nE_{j-1} \\
&\qquad \times (\mathbf{x}_n^*\breve{A}_{ij}^{-1}(z_2)T_nA_{ij}^{-1}(z_1)\mathbf{x}_n\operatorname{tr}T_n\breve{A}_{ij}^{-1}(z_2)T_nT^{-1}(z_1)) + o_p(1).
\end{aligned}
$$

Reversing the above procedure, one finds that we may also replace $A_{ij}^{-1}(z_1)$ and $\breve{A}_{ij}^{-1}(z_2)$ in $B(z_1,z_2)$ by $A_j^{-1}(z_1)$ and $A_j^{-1}(z_2)$, respectively. That is,

$$
\begin{aligned}
B(z_1,z_2) &= -\frac{b_{n1}(z_1)b_{n1}(z_2)(j-1)}{N^2}E_{j-1}\mathbf{x}_n^*A_j^{-1}(z_1)T_n\breve{A}_j^{-1}(z_2)\mathbf{x}_nE_{j-1} \\
&\quad \times (\mathbf{x}_n^*\breve{A}_j^{-1}(z_2)T_nA_j^{-1}(z_1)\mathbf{x}_n\operatorname{tr}T_n\breve{A}_j^{-1}(z_2)T_nT^{-1}(z_1)) + o_p(1).
\end{aligned}
$$

Using the martingale decomposition (4.8), one can further show that

$$B(z_1,z_2) = -\frac{b_{n1}(z_1)b_{n1}(z_2)(j-1)}{N^2} E_{j-1}\mathbf{x}_n^* A_j^{-1}(z_1)T_n\breve{A}_j^{-1}(z_2)\mathbf{x}_n \times E_{j-1}\mathbf{x}_n^*\breve{A}_j^{-1}(z_2)T_nA_j^{-1}(z_1)\mathbf{x}_n E_{j-1}\operatorname{tr} T_n\breve{A}_j^{-1}(z_2)T_nT^{-1}(z_1) + o_p(1).$$

It is easy to verify that

$$N^{-1}\operatorname{tr}(T_nM(z_2)T_nT^{-1}(z_1)) = o_p(1)$$

when $M(z_2)$ takes the value $\breve{B}_j(z_2)$, $\breve{C}_j(z_2)$ or $\breve{D}_j(z_2)$. Thus, substituting the decomposition (4.10) for $\breve{A}_j^{-1}(z_2)$ in the above approximation for $B(z_1,z_2)$, one finds that

$$B(z_1,z_2) = \frac{b_{n1}(z_1)b_{n1}(z_2)(j-1)}{N^2} E_{j-1}\mathbf{x}_n^* A_j^{-1}(z_1)T_n\breve{A}_j^{-1}(z_2)\mathbf{x}_n \times E_{j-1}\mathbf{x}_n^*\breve{A}_j^{-1}(z_2)T_nA_j^{-1}(z_1)\mathbf{x}_n \times E_{j-1}\operatorname{tr} T_nT^{-1}(z_2)T_nT^{-1}(z_1) + o_p(1). \tag{4.18}$$

Finally, let us consider the first term of (4.11). Using the expression for $A_j^{-1}(z_1)$ in (4.10), we obtain

$$\begin{aligned} &-E_{j-1}\mathbf{x}_n^*A_j^{-1}(z_1)T_n\breve{A}_j^{-1}(z_2)\mathbf{x}_nE_{j-1}\mathbf{x}_n^*\breve{A}_j^{-1}(z_2)T_nT^{-1}(z_1)\mathbf{x}_n \\ &\quad = W_1(z_1,z_2) + W_2(z_1,z_2) + W_3(z_1,z_2) + W_4(z_1,z_2), \end{aligned} \tag{4.19}$$

where

$$W_1(z_1,z_2) = E_{j-1}\mathbf{x}_n^*T^{-1}(z_1)T_n\breve{A}_j^{-1}(z_2)\mathbf{x}_nE_{j-1}\mathbf{x}_n^*\breve{A}_j^{-1}(z_2)T_nT^{-1}(z_1)\mathbf{x}_n,$$

$$W_2(z_1,z_2) = -b_{n1}(z_1)E_{j-1}\mathbf{x}_n^*B_{j2}^{-1}(z_1)T_n\breve{A}_j^{-1}(z_2)\mathbf{x}_n \times E_{j-1}\mathbf{x}_n^*\breve{A}_j^{-1}(z_2)T_nT^{-1}(z_1)\mathbf{x}_n,$$

$$W_3(z_1,z_2) = -E_{j-1}\mathbf{x}_n^*C_j^{-1}(z_1)T_n\breve{A}_j^{-1}(z_2)\mathbf{x}_nE_{j-1}\mathbf{x}_n^*\breve{A}_j^{-1}(z_2)T_nT^{-1}(z_1)\mathbf{x}_n$$

and

$$W_4(z_1,z_2) = -E_{j-1}\mathbf{x}_n^*D_j^{-1}(z_1)T_n\breve{A}_j^{-1}(z_2)\mathbf{x}_nE_{j-1}\mathbf{x}_n^*\breve{A}_j^{-1}(z_2)T_nT^{-1}(z_1)\mathbf{x}_n.$$

By the same argument as (4.12), one can obtain

$$E|W_3(z_1,z_2)| = o(1) \quad \text{and} \quad E|W_4(z_1,z_2)| = o(1). \tag{4.20}$$

Furthermore, as in dealing with $B(z_1,z_2)$, the first $\breve{A}_j^{-1}(z_2)$ in $W_2(z_1,z_2)$ can be replaced by $-b_{n1}(z_2)\breve{A}_{ij}^{-1}(z_2)\mathbf{s}_i\mathbf{s}_i^*\breve{A}_{ij}^{-1}(z_2)$, that is,

$W_2(z_1,z_2)$

$$
\begin{aligned}
&= b_{n1}(z_1)b_{n1}(z_2) \\
&\quad \times \sum_{i<j} E_{j-1}\mathbf{x}_n^* T^{-1}(z_1)(\mathbf{s}_i\mathbf{s}_i^* - N^{-1}T_n)A_{ij}^{-1}(z_1)T_n \\
&\qquad \times \breve{A}_{ij}^{-1}(z_2)\mathbf{s}_i\mathbf{s}_i^*\breve{A}_{ij}^{-1}(z_2)\mathbf{x}_n E_{j-1}\mathbf{x}_n^*\breve{A}_j^{-1}(z_2)T_nT^{-1}(z_1)\mathbf{x}_n + o_p(1) \\
&= b_{n1}(z_1)b_{n1}(z_2) \\
&\quad \times \sum_{i<j} E_{j-1}\mathbf{x}_n^* T^{-1}(z_1)\mathbf{s}_i\mathbf{s}_i^* A_{ij}^{-1}(z_1)T_n \\
&\qquad \times \breve{A}_{ij}^{-1}(z_2)\mathbf{s}_i\mathbf{s}_i^*\breve{A}_{ij}^{-1}(z_2)\mathbf{x}_n E_{j-1}\mathbf{x}_n^*\breve{A}_j^{-1}(z_2)T_nT^{-1}(z_1)\mathbf{x}_n + o_p(1) \\
&= \frac{b_{n1}(z_1)b_{n1}(z_2)(j-1)}{N^2} \\
&\quad \times E_{j-1}(\mathbf{x}_n^*T^{-1}(z_1)T_n\breve{A}_j^{-1}(z_2)\mathbf{x}_n \operatorname{tr} T_nA_j^{-1}(z_1)T_n\breve{A}_j^{-1}(z_2)) \\
&\quad \times E_{j-1}\mathbf{x}_n^*\breve{A}_j^{-1}(z_2)T_nT^{-1}(z_1)\mathbf{x}_n + o_p(1).
\end{aligned}
$$

It can also be verified that

$$\mathbf{x}_n^* M(z_2)T_nT^{-1}(z_1)\mathbf{x}_n = o_p(1),$$

when $M(z_2)$ takes the value $\breve{B}_j(z_2)$, $\breve{C}_j(z_2)$ or $\breve{D}_j(z_2)$. Therefore, $W_2(z_1, z_2)$ can be further approximated by

$$
\begin{aligned}
W_2(z_1,z_2) &= \frac{b_{n1}(z_1)b_{n1}(z_2)(j-1)}{N^2}(\mathbf{x}_n^*T^{-1}(z_1)T_nT^{-1}(z_2)\mathbf{x}_n)^2 \\
&\quad \times E_{j-1}\operatorname{tr}(T_nA_j^{-1}(z_1)T_n\breve{A}^{-1}(z_2)) + o_p(1).
\end{aligned} \tag{4.21}
$$

In (2.18) of [7], it is proved that

$$
\begin{aligned}
&E_{j-1}\operatorname{tr}(T_nA_j^{-1}(z_1)T_n\breve{A}_j^{-1}(z_2)) \\
&\quad = \frac{\operatorname{tr}(T_nT^{-1}(z_1)T_nT^{-1}(z_2)) + o_p(1)}{1-(j-1)/N^2 z_1z_2\underline{m}(z_1)\underline{m}(z_2)\operatorname{tr}(T_nT^{-1}(z_1)T_nT^{-1}(z_2))}.
\end{aligned}
$$

By the same method, $W_1(z_1, z_2)$ can be approximated by

$$W_1(z_1,z_2) = \mathbf{x}_n^*T^{-1}(z_1)T_nT^{-1}(z_2)\mathbf{x}_n\mathbf{x}_n^*T^{-1}(z_2)T_nT^{-1}(z_1)\mathbf{x}_n + o_p(1). \tag{4.22}$$

Consequently, from (4.11)–(4.22), we obtain

$$
\begin{aligned}
&E_{j-1}\mathbf{x}_n^*A_j^{-1}(z_1)T_n\breve{E}_j(A_j^{-1}(z_2)\mathbf{x}_nE_{j-1}\mathbf{x}_n^*\breve{A}_j^{-1}(z_2))T_nA_j^{-1}(z_1)\mathbf{x}_n \\
&\quad \times\left[1-\frac{j-1}{N}b_{n1}(z_1)b_{n1}(z_2)\frac{1}{N}\operatorname{tr}T^{-1}(z_2)T_nT^{-1}(z_1)T_n\right] \\
&\qquad = \mathbf{x}_n^*T^{-1}(z_1)T_nT^{-1}(z_2)\mathbf{x}_n\mathbf{x}_n^*T^{-1}(z_2)T_nT^{-1}(z_1)\mathbf{x}_n
\end{aligned} \tag{4.23}
$$

$$\times\Big(1+\frac{j-1}{N}b_{n1}(z_1)b_{n1}(z_2)\frac{1}{N}E_{j-1}\operatorname{tr}(A_j^{-1}(z_1)T_n\breve{A}_j^{-1}(z_2)T_n)\Big)$$
$$+ o_p(1).$$

Recall that $b_{n1}(z)\to -z\underline{m}(z)$ and $F^{T_n}\to H$. Hence,

$$\begin{aligned}
d(z_1,z_2) &:= \lim b_{n1}(z_1)b_{n1}(z_2)\frac{1}{N}\operatorname{tr}(T^{-1}(z_1)T_nT^{-1}(z_2)T_n)\\
&= \int \frac{ct^2\underline{m}(z_1)\underline{m}(z_2)}{(1+t\underline{m}(z_1))(1+t\underline{m}(z_2))}\,dH(t) \qquad (4.24)\\
&= 1+\frac{\underline{m}(z_1)\underline{m}(z_2)(z_1-z_2)}{\underline{m}(z_2)-\underline{m}(z_1)}.
\end{aligned}$$

By the conditions of Theorem 2,

$$\begin{aligned}
h(z_1,z_2) &= \lim z_1z_2\underline{m}(z_1)\underline{m}(z_2)\mathbf{x}_n^*T^{-1}(z_1)T_nT^{-1}(z_2)\mathbf{x}_n\\
&\qquad \times \mathbf{x}_n^*T^{-1}(z_2)T_nT^{-1}(z_1)\mathbf{x}_n\\
&= \frac{\underline{m}(z_1)\underline{m}(z_2)}{z_1z_2}\bigg(\int\frac{t^2\underline{m}(z_1)\underline{m}(z_2)}{(1+t\underline{m}(z_1))(1+t\underline{m}(z_2))}\,dH(t)\bigg)^2 \qquad (4.25)\\
&= \frac{\underline{m}(z_1)\underline{m}(z_2)}{z_1z_2}\bigg(\int\frac{t\,dH(t)}{(1+t\underline{m}(z_1))(1+t\underline{m}(z_2))}\bigg)^2\\
&= \frac{\underline{m}(z_1)\underline{m}(z_2)}{z_1z_2}\bigg(\frac{z_1m(z_1)-z_2m(z_2)}{(\underline{m}(z_2)-\underline{m}(z_1))}\bigg)^2.
\end{aligned}$$

From (4.10), (4.24) and (4.25), we get

$$\begin{aligned}
(4.7) &\xrightarrow{i.p.} h(z_1,z_2)\bigg(\int_0^1\frac{1}{(1-td(z_1,z_2))}\,dt+\int_0^1\frac{td(z_1,z_2)}{(1-td(z_1,z_2))^2}\,dt\bigg)\\
&= \frac{h(z_1,z_2)}{1-d(z_1,z_2)} = \frac{(z_2\underline{m}(z_2)-z_1\underline{m}(z_1))^2}{c^2z_1z_2(z_2-z_1)(\underline{m}(z_2)-\underline{m}(z_1))}.
\end{aligned}$$

**5. Tightness of $M_n^1(z)$ and convergence of $M_n^2(z)$.** First, we proceed to the proof of the tightness of $M_n^1(z)$. By (4.3), we obtain

$$E\bigg|\sum_{i=1}^r a_i\sum_{j=1}^N Y_j(z_i)\bigg|^2=\sum_{j=1}^N E\bigg|\sum_{i=1}^r a_iY_j(z_i)\bigg|^2\le K.$$

Thus, as pointed out in [7], condition (i) of Theorem 12.3 of [8] is satisfied. Therefore, to complete the proof of tightness, we only need verify that

$$E\frac{|M_n^1(z_1)-M_n^1(z_2)|^2}{|z_1-z_2|^2}\le K \qquad \text{if } z_1,z_2\in\mathcal{C}. \tag{5.1}$$

Write

$$Q(z_1,z_2)=\sqrt{N}\mathbf{x}_n(A_n-z_1I)^{-1}(A_n-z_2I)^{-1}\mathbf{x}_n.$$

Recalling the definition of $M_n^1$, we have

$$\begin{aligned}&\frac{|M_n^1(z_1)-M_n^1(z_2)|}{|z_1-z_2|}\\&\qquad\le\begin{cases}|Q(z_1,z_2)-EQ(z_1,z_2)|,\\\qquad\text{if } z_1,z_2\in\mathcal{C}_0\cup\mathcal{C}_0,\\|Q(z_1,z_{2+})-EQ(z_1,z_{2+})|+|Q(z_1,z_{2-})-EQ(z_1,z_{2-})|,\\\qquad\text{if } z_1\in\mathcal{C}_0\cup\mathcal{C}_0 \text{ only},\\|Q(z_+,z_-)-EQ(z_+,z_-)|,\\\qquad\text{otherwise},\end{cases}\end{aligned}$$

where $\Re(z_{2\pm})=\Re z_2$, $\Im(z_{2\pm})=\pm\delta_n n^{-1}$, $\Im(z_\pm)=\pm\delta_n n^{-1}$ and $\Re(z_\pm)=u_l$ or $u_r$. Thus we need only to show (5.1) when $z_1,z_2\in\mathcal{C}_0\cup\bar{\mathcal{C}}_0$.

From the identity above (3.7) in [7], we obtain

$$\begin{aligned}\frac{M_n^1(z_1)-M_n^1(z_2)}{z_1-z_2}&=\sqrt{N}\sum_{j=1}^N(E_j-E_{j-1})\mathbf{x}_n^*A^{-1}(z_1)A^{-1}(z_2)\mathbf{x}_n\\&=V_1(z_1,z_2)+V_2(z_1,z_2)+V_3(z_1,z_2),\end{aligned}\tag{5.2}$$

where

$$\begin{aligned}V_1(z_1,z_2)=\sqrt{N}\sum_{j=1}^N&(E_j-E_{j-1})\beta_j(z_1)\beta_j(z_2)\mathbf{s}_j^*A_j^{-1}(z_1)A_j^{-1}(z_2)\mathbf{s}_j\mathbf{s}_j^*\\&\times A_j^{-1}(z_2)\mathbf{x}_n\mathbf{x}_n^*A_j^{-1}(z_1)\mathbf{s}_j,\end{aligned}$$

$$V_2(z_1,z_2)=\sqrt{N}\sum_{j=1}^N(E_j-E_{j-1})\beta_j(z_1)\mathbf{s}_j^*A_j^{-1}(z_1)A_j^{-1}(z_2)\mathbf{x}_n\mathbf{x}_n^*A_j^{-1}(z_1)\mathbf{s}_j$$

and

$$V_3(z_1,z_2)=\sqrt{N}\sum_{j=1}^N(E_j-E_{j-1})\beta_j(z_2)\mathbf{s}_j^*A_j^{-1}(z_2)\mathbf{x}_n\mathbf{x}_n^*A_j^{-1}(z_1)A_j^{-1}(z_2)\mathbf{s}_j.$$

Applying (3.1) and the bounds for $\beta_j(z)$ and $\mathbf{s}_j^*A_j^{-1}(z_1)A_j^{-1}(z_2)\mathbf{s}_j$ given in the remark concerning (3.2) in [7], we obtain

$$\begin{aligned}E|V_1(z_1,z_2)|^2=N\sum_{j=1}^N&E|(E_j-E_{j-1})\beta_j(z_1)\beta_j(z_2)\mathbf{s}_j^*A_j^{-1}(z_1)\\&\times A_j^{-1}(z_2)\mathbf{s}_j\mathbf{s}_j^*A_j^{-1}(z_2)\mathbf{x}_n\mathbf{x}_n^*\times A_j^{-1}(z_1)\mathbf{s}_j|^2\end{aligned}$$

$$\begin{aligned}(5.3)\qquad &\le KN^2(E|\mathbf{s}_j^* A_j^{-1}(z_2)\mathbf{x}_n\mathbf{x}_n^* A_j^{-1}(z_1)\mathbf{s}_j|^2 \\ &\qquad + v^{-12}n^2 P(\|A_1\| > u_r \text{ or } \lambda_{\min}^{A_{(1)}} < u_l)) \\ &\le K,\end{aligned}$$

where (1.9a) and (1.9b) in [7] are also employed. It is easy to see that (1.9a) and (1.9b) in [7] also hold under our truncation case. Similarly, the above argument can also be used to treat $V_2(z_1, z_2)$ and $V_3(z_1, z_2)$. Therefore, we have completed the proof of (5.1).

Next, we will consider the convergence of $M_n^2(z)$. Note that

$$m_{F^{c_n,H_n}}(z) = -\frac{1}{z}\int \frac{1}{1 + t\underline{m}_{F^{c_n,H_n}}(z)}\, dH_n(t). \tag{5.4}$$

Substituting (2.5) into (5.4), we obtain that

$$\begin{aligned}&z\sqrt{N}((\mathbf{x}_n^* E(A^{-1})(z)\mathbf{x}_n - m_{F^{c_n,H_n}}(z))) \\ (5.5)\qquad &= \sqrt{N}\left(\mathbf{x}_n^*(-E\underline{m}_n(z)T_n - I)^{-1}\mathbf{x}_n + \int \frac{1}{1 + t\underline{m}_{F^{c_n,H_n}}(z)}\, dH_n(t)\right) \\ &\quad + \sqrt{N}z(\delta_1 + \delta_2 + \delta_3).\end{aligned}$$

Applying (2.6)–(2.8), we have

$$\sqrt{N}z(\delta_1 + \delta_2 + \delta_3) = o(1).$$

On the other hand, in Section 4 of [7], it is proved that

$$\sup_z \sqrt{N}(\underline{m}_{F^{c_n,H_n}}(z) - E\underline{m}_n(z)) \to 0. \tag{5.6}$$

Following a similar line to (4.3) of [7], along with (4.2) of [7], we can obtain

$$\sup_{n, z\in\mathcal{C}_0} \|(\underline{m}_{F^{c_n,H_n}}(z)T_n + I)^{-1}\| < \infty.$$

It follows, via (4.3) of [7] and the assumption of Theorem 2, that

$$\sup_{n, z\in\mathcal{C}_0} \int \frac{t}{(1 + t\underline{m}_{F^{c_n,H_n}}(z))(tE\underline{m}_n(z) + 1)}\, dH_n(t) < \infty. \tag{5.7}$$

Appealing to (4.1), (4.3) in [7], (5.6) and (5.7), we conclude that

$$\begin{aligned}&\sqrt{N}\left(\mathbf{x}_n^*(E\underline{m}_n(z)T_n + I)^{-1}\mathbf{x}_n - \int \frac{1}{E\underline{m}_n(z)t + 1}\, dH_n(t)\right) \\ &\quad = \sqrt{N}\left(\mathbf{x}_n^*(\underline{m}_{F^{c_n,H_n}}(z)T_n + I)^{-1}\mathbf{x}_n - \int \frac{1}{\underline{m}_{F^{c_n,H_n}}(z)t + 1}\, dH_n(t)\right) + o(1).\end{aligned}$$

Using (5.6) and (5.7), we also have

$$\begin{aligned}
&\sqrt{N}\left(\int \frac{1}{E\underline{m}_n(z)t+1}\,dH_n(t) - \int \frac{1}{1+t\underline{m}_{F^{c_n,H_n}}(z)}\,dH_n(t)\right)\\
&\quad = \sqrt{N}(\underline{m}_{F^{c_n,H_n}}(z) - E\underline{m}_n(z)) \int \frac{t}{(1+t\underline{m}_{F^{c_n,H_n}}(z))(tE\underline{m}_n(z)+1)}\,dH_n(t)\\
&\quad = o(1).
\end{aligned}$$

Combining the above arguments, we can conclude that

$$(5.5) \to 0.$$

**6. Proof of Theorem 3 and supplement to Remark 7.** In this section, when $T_n = I$, we will show that formula (1.12) includes (1.2) of [15] as a special case and we will also present a proof of (1.14).

First, one can easily verify that (1.15) reduces to (1.2) of [15] when $g(x) = x^r$. Next, our goal is to prove that (1.12) implies (1.14) under the condition of Theorem 3.

Write

$$\begin{aligned}
(z_2\underline{m}(z_2) - z_1\underline{m}(z_1))^2 &= z_1z_2(\underline{m}(z_1) - \underline{m}(z_2))^2 + \underline{m}(z_1)\underline{m}(z_2)(z_2 - z_1)^2\\
&\quad + z_2\underline{m}(z_2)(z_2 - z_1)(\underline{m}(z_2) - \underline{m}(z_1))\\
&\quad + z_1\underline{m}(z_1)(z_2 - z_1)(\underline{m}(z_2) - \underline{m}(z_1)).
\end{aligned} \tag{6.1}$$

Recall that

$$z = -\frac{1}{\underline{m}(z)} + c\int \frac{t}{1+t\underline{m}(z)}\,dH(t), \tag{6.2}$$

from which [together with assumption (1.13)] we obtain

$$\begin{aligned}
\frac{\underline{m}(z_1)\underline{m}(z_2)(z_2 - z_1)}{\underline{m}(z_2) - \underline{m}(z_1)} &= 1 - c\int \frac{t^2\underline{m}(z_1)\underline{m}(z_2)\,dH(t)}{(1+t\underline{m}(z_1))(1+t\underline{m}(z_2))}\\
&= 1 - c\int \frac{t\underline{m}(z_1)}{1+t\underline{m}(z_1)}\,dH(t) \int \frac{t\underline{m}(z_2)}{1+t\underline{m}(z_2)}\,dH(t)\\
&= 1 - c^{-1}(1 + z_1\underline{m}(z_1))(1 + z_2\underline{m}(z_2)).
\end{aligned}$$

Replacing one copy of $z_2 - z_1$ by this in the second term on the right-hand side of (6.1), we obtain

$$\begin{aligned}
&(z_2\underline{m}(z_2) - z_1\underline{m}(z_1))^2\\
&\quad = z_1z_2(\underline{m}(z_1) - \underline{m}(z_2))^2\\
&\qquad + (z_2 - z_1)(\underline{m}(z_2) - \underline{m}(z_1))[(1 - c^{-1})(1 + z_1\underline{m}(z_1) + z_2\underline{m}(z_2))\\
&\qquad\qquad - c^{-1}z_1z_2\underline{m}(z_1)\underline{m}(z_2)].
\end{aligned} \tag{6.3}$$

Using this and the facts that (1), (2), (3), we obtain RHS

$$\frac{1}{2\pi i}\oint_{\mathcal{C}} z^{-1}g(z)\,dz = \begin{cases} g(0), & \text{if } \mathcal{C} \text{ encloses the origin,} \\ 0, & \text{otherwise,} \end{cases}$$

$$\frac{1}{2\pi i}\oint_{\mathcal{C}} \underline{m}(z)g(z)\,dz = -\int g(x)\,d\underline{F}^{c,H}(x),$$

$$-\frac{1}{4\pi^2}\oint_{\mathcal{C}_1}\oint_{\mathcal{C}_2} \frac{\underline{m}(z_2)-\underline{m}(z_1)}{z_2-z_1} g_1(z_1)g_2(z_2)\,dz_1\,dz_2 = \int g_1(x)g_2(x)\,d\underline{F}^{c,H}(x),$$

when $\mathcal{C}_j$ enclose the origin, we obtain

RHS of (1.12)

$$\begin{aligned}
&= \frac{2}{c^2}\int g_1(x)g_2(x)\,d\underline{F}^{c,H}(x) \\
&\quad + \frac{2(c-1)}{c^3}\biggl(g_1(0)g_2(0) - g_1(0)\int g_2(x)\,d\underline{F}^{c,H}(x) \\
&\qquad\qquad - g_2(0)\int g_1(x)\,d\underline{F}^{c,H}(x)\biggr) \\
&\quad - \frac{2}{c^3}\int g_1(x)\,d\underline{F}^{c,H}(x)\int g_2(x)\,d\underline{F}^{c,H}(x) \\
&= \frac{2}{c}\biggl(\int g_1(x)g_2(x)\,dF^{c,H}(x) - \int g_1(x)\,dF^{c,H}(x)\int g_2(x)\,dF^{c,H}(x)\biggr).
\end{aligned}$$

The same we can obtain when $\mathcal{C}_j$ does not enclose the origin, we can obtain the same result even more easily.

Thus the result is proved.

## 7. Truncation, recentralization and renormalization of $X_{ij}$.

In this section, to facilitate the analysis in the other sections, the underlying random variables need to be truncated, recentralized and renormalized. Since the argument for (1.8) in [7] can be carried directly over to the present case, we can then select $\varepsilon_n$ such that

$$\varepsilon_n \to 0 \quad \text{and} \quad \varepsilon_n^{-4}\int_{\{|X_{11}|\ge \varepsilon_n n^{1/4}\}} |X_{11}|^4 \to 0. \tag{7.1}$$

Let $\hat{X}_{ij} = X_{ij}I(|X_{ij}| \le \varepsilon_n n^{1/4}) - EX_{ij}I(|X_{ij}| \le \varepsilon_n n^{1/4})$, $\tilde{X}_n = X_n - \hat{X}_n$, and $\hat{A}_n = \frac{1}{N}T_n^{1/2}\hat{X}_n\hat{X}_n^* T_n^{1/2}$, where $\hat{X}_n = (\hat{X}_{ij})$. Let $\sigma_n^2 = E|\hat{X}_{11}|^2$ and $\check{A}_n = \frac{1}{N\sigma_n^2}T_n^{1/2}\hat{X}_n\hat{X}_n^*T_n^{1/2}$. Write $A^{-1}(z) = (A_n - zI)^{-1}$ and $\check{A}^{-1}(z) = (\check{A}_n - zI)^{-1}$.

LEMMA 4. *Suppose that $X_{ij} \in \mathcal{C}, i = 1, \ldots, n, j = 1, \ldots, N$, are i.i.d. with $EX_{11} = 0, E|X_{11}|^2 = 1$ and $E|X_{11}|^4 < \infty$. For $z \in \mathcal{C}_0$, we then have*

$$\sqrt{N}\mathbf{x}_n^*(\check{A}^{-1}(z)\mathbf{x}_n - \mathbf{x}_n^* A^{-1}(z))\mathbf{x}_n \xrightarrow{i.p.} 0. \tag{7.2}$$

PROOF. Corresponding to the truncated and renormalized variables, we similarly define $\hat{\mathbf{s}}_j$, $\tilde{\mathbf{s}}_j$, $\hat{A}^{-1}(z)$, $\hat{A}_j^{-1}(z)$, $\hat{\beta}_j(z)$, $\hat{\beta}_{j_1j_2}(z)$ and $\hat{\beta}_{j_2j_1}(z)$.

We then have

$$\begin{aligned}
&\sqrt{N}(\mathbf{x}_n^*\hat{A}^{-1}(z)\mathbf{x}_n - \mathbf{x}_n^* A^{-1}(z)\mathbf{x}_n)\\
&\quad= \frac{1}{\sqrt{N}}(\mathbf{x}_n^*\hat{A}^{-1}(z)T_n^{1/2}(X_n - \hat{X}_n)X_n^* T_n^{1/2}A^{-1}(z)\mathbf{x}_n\\
&\qquad\qquad + \mathbf{x}_n^*\hat{A}^{-1}(z)T_n^{1/2}\hat{X}_n(X_n^* - \hat{X}_n^*)T_n^{1/2}A^{-1}(z)\mathbf{x}_n)\\
&\quad= \sqrt{N}\sum_{j=1}^N \mathbf{x}_n^*\hat{A}^{-1}(z)\tilde{\mathbf{s}}_j\mathbf{s}_j^* A^{-1}(z)\mathbf{x}_n + \sqrt{N}\sum_{j=1}^N \mathbf{x}_n^*\hat{A}^{-1}(z)\hat{\mathbf{s}}_j\tilde{\mathbf{s}}_j^* A^{-1}(z)\mathbf{x}_n \qquad (7.3)\\
&\quad= \sqrt{N}\sum_{j=1}^N \beta_j(z)\mathbf{x}_n^*\hat{A}^{-1}(z)\tilde{\mathbf{s}}_j\mathbf{s}_j^* A_j^{-1}(z)\mathbf{x}_n\\
&\qquad + \sqrt{N}\sum_{j=1}^N \hat{\beta}_j(z)\mathbf{x}_n^*\hat{A}_j^{-1}(z)\hat{\mathbf{s}}_j\tilde{\mathbf{s}}_j^* A^{-1}(z)\mathbf{x}_n\\
&\quad\stackrel{\triangle}{=} \omega_1 + \omega_2.
\end{aligned}$$

Consider first the term $\omega_1$.

$$\begin{aligned}
E|\omega_1|^2 &\le N\sum_{j=1}^N E|\beta_j(z)\mathbf{x}_n^*\hat{A}^{-1}(z)\tilde{\mathbf{s}}_j\mathbf{s}_j^* A_j^{-1}(z)\mathbf{x}_n|^2\\
&\qquad + N\sum_{j_1\ne j_2}|E\beta_{j_1}(z)\mathbf{x}_n^*\hat{A}^{-1}(z)\tilde{\mathbf{s}}_{j_1}\mathbf{s}_{j_1}^* A_{j_1}^{-1}(z)\mathbf{x}_n \\
&\qquad\qquad \times \beta_{j_2}(\bar{z})\mathbf{x}_n^* A_{j_2}^{-1}(\bar{z})\mathbf{s}_{j_2}\tilde{\mathbf{s}}_{j_2}^*\hat{A}^{-1}(\bar{z})\mathbf{x}_n| \qquad (7.4)\\
&\stackrel{\triangle}{=} \omega_{11} + \omega_{12},
\end{aligned}$$

where $\beta_{j_2}(\bar{z})$ is the complex conjugate of $\beta_{j_2}(z)$.

Our next goal is to show that the above two terms converge to zero for all $z \in \mathcal{C}_u$ and $z \in \mathcal{C}_l$ when $u_l < 0$.

In this case, $\beta_i(z)$ is bounded. It is straightforward to verify that

$$\begin{aligned}
&E|\beta_j(z)\mathbf{x}_n^*\hat{A}_j^{-1}(z)\tilde{\mathbf{s}}_j\mathbf{s}_j^* A_j^{-1}(z)\mathbf{x}_n|^2\\
&\qquad \le KE|\mathbf{s}_j^* A_j^{-1}(z)\mathbf{x}_n\mathbf{x}_n^*\hat{A}_j^{-1}(z)\tilde{\mathbf{s}}_j|^2 = o(N^{-2}),
\end{aligned}$$

$$E|\hat{\mathbf{s}}_j^* \hat{A}_j^{-1}(z)\tilde{\mathbf{s}}_j - E(\hat{X}_{11}\tilde{X}_{11}^*)N^{-1}\operatorname{tr}\hat{A}_j^{-1}(z)T_n|^4 = o(N^{-2}),$$
$$E|\mathbf{s}_j^* A_j^{-1}(z)\mathbf{x}_n\mathbf{x}_n^*\hat{A}_j^{-1}(z)\hat{\mathbf{s}}_j|^4 = O(N^{-3}). \tag{7.5}$$

Therefore,

$$\begin{aligned} &E|\beta_j(z)\hat{\beta}_j(z)\mathbf{x}_n^*\hat{A}_j^{-1}(z)\hat{\mathbf{s}}_j\hat{\mathbf{s}}_j^*\hat{A}_j^{-1}(z)\tilde{\mathbf{s}}_j\mathbf{s}_j^*A_j^{-1}(z)\mathbf{x}_n|^2 \\ &\quad \le K(E|\hat{\mathbf{s}}_j^*\hat{A}_j^{-1}(z)\tilde{\mathbf{s}}_j - E(\hat{X}_{11}\tilde{X}_{11}^*)N^{-1}\operatorname{tr}\hat{A}_j^{-1}(z)T_n|^4)^{1/2} \\ &\qquad \times (E|\mathbf{s}_j^*A_j^{-1}(z)\mathbf{x}_n\mathbf{x}_n^*\hat{A}_j^{-1}(z)\hat{\mathbf{s}}_j|^4)^{1/2} \\ &\qquad + E|X_{11}|^2 I(|X_{11}| \ge \varepsilon_n n^{1/4})E|\mathbf{s}_j^*A_j^{-1}(z)\mathbf{x}_n\mathbf{x}_n^*\hat{A}_j^{-1}(z)\hat{\mathbf{s}}_j|^2 = o(N^{-2}). \end{aligned}$$

It then follows that

$$\begin{aligned} \omega_{11} &\le KN\sum_{j=1}^{N}(E|\beta_j(z)\hat{\beta}_j(z)\mathbf{x}_n^*\hat{A}_j^{-1}(z)\hat{\mathbf{s}}_j\hat{\mathbf{s}}_j^*\hat{A}_j^{-1}(z)\tilde{\mathbf{s}}_j\mathbf{s}_j^*A_j^{-1}(z)\mathbf{x}_n|^2 \\ &\qquad\qquad + E|\beta_j(z)\mathbf{x}_n^*\hat{A}_j^{-1}(z)\tilde{\mathbf{s}}_j\mathbf{s}_j^*A_j^{-1}(z)\mathbf{x}_n|^2) \qquad (7.6) \\ &= o(1). \end{aligned}$$

Now consider the term $\omega_{12}$. We have

$$\begin{aligned} &E\beta_{j_1}(z)\beta_{j_2}(\bar{z})\mathbf{s}_{j_1}^*A_{j_1}^{-1}(z)\mathbf{x}_n\mathbf{x}_n^*A_{j_2}^{-1}(\bar{z})\mathbf{s}_{j_2}\tilde{\mathbf{s}}_{j_2}^*\hat{A}^{-1}(\bar{z})\mathbf{x}_n\mathbf{x}_n^*\hat{A}^{-1}(z)\tilde{\mathbf{s}}_{j_1} \\ &\quad = \Delta_1 + \Delta_2 + \Delta_3 + \Delta_4, \end{aligned}$$

where

$$\begin{aligned} \Delta_1 &= E\beta_{j_1}(z)\beta_{j_2}(\bar{z})\hat{\beta}_{j_1}(z)\hat{\beta}_{j_2}(\bar{z})\mathbf{s}_{j_1}^*A_{j_1}^{-1}(z)\mathbf{x}_n\mathbf{x}_n^*A_{j_2}^{-1}(\bar{z})\mathbf{s}_{j_2}\tilde{\mathbf{s}}_{j_2}^*\hat{A}_{j_2}^{-1}(\bar{z})\hat{\mathbf{s}}_{j_2} \\ &\quad \times \hat{\mathbf{s}}_{j_2}^*\hat{A}_{j_2}^{-1}(\bar{z})\mathbf{x}_n\mathbf{x}_n^*\hat{A}_{j_1}^{-1}(z)\hat{\mathbf{s}}_{j_1}\hat{\mathbf{s}}_{j_1}^*\hat{A}_{j_1}^{-1}(z)\tilde{\mathbf{s}}_{j_1}, \\ \Delta_2 &= -E\beta_{j_1}(z)\beta_{j_2}(\bar{z})\hat{\beta}_{j_2}(z)\mathbf{s}_{j_1}^*A_{j_1}^{-1}(z)\mathbf{x}_n\mathbf{x}_n^*A_{j_2}^{-1}(\bar{z})\mathbf{s}_{j_2} \\ &\quad \times \tilde{\mathbf{s}}_{j_2}^*\hat{A}_{j_2}^{-1}(\bar{z})\hat{\mathbf{s}}_{j_2}\hat{\mathbf{s}}_{j_2}^*\hat{A}_{j_2}^{-1}(\bar{z})\mathbf{x}_n\mathbf{x}_n^*\hat{A}_{j_1}^{-1}(z)\tilde{\mathbf{s}}_{j_1}, \\ \Delta_3 &= -E\beta_{j_1}(z)\beta_{j_2}(\bar{z})\hat{\beta}_{j_1}(z)s_{j_1}^*A_{j_1}^{-1}(z)\mathbf{x}_n\mathbf{x}_n^*A_{j_2}^{-1}(\bar{z})\mathbf{s}_{j_2}\tilde{\mathbf{s}}_{j_2}^*\hat{A}_{j_2}^{-1}(\bar{z})\mathbf{x}_n \\ &\quad \times \mathbf{x}_n^*\hat{A}_{j_1}^{-1}(z)\hat{\mathbf{s}}_{j_1}\hat{\mathbf{s}}_{j_1}^*\hat{A}_{j_1}^{-1}(z)\tilde{\mathbf{s}}_{j_1}, \\ \Delta_4 &= E\beta_{j_1}(z)\beta_{j_2}(\bar{z})\mathbf{s}_{j_1}^*A_{j_1}^{-1}(z)\mathbf{x}_n\mathbf{x}_n^*A_{j_2}^{-1}(\bar{z})\mathbf{s}_{j_2}\tilde{\mathbf{s}}_{j_2}^*\hat{A}_{j_2}^{-1}(\bar{z})\mathbf{x}_n\mathbf{x}_n^*\hat{A}_{j_1}^{-1}(z)\tilde{\mathbf{s}}_{j_1}. \end{aligned}$$

In the sequel, $\Delta_1$ will be further decomposed. However, since the expansions of $\Delta_1$ are rather complicated, as an illustration, we only present estimates for some typical terms; other terms can be estimated similarly. The main technique is to extract the $j_1$th and $j_2$th column, respectively, from $T_n^{1/2}\hat{X}_n$ of $A_{j_2}^{-1}(z)$ and $A_{j_1}^{-1}(z)$, and so on. We will evaluate the following terms of $\Delta_1$:

$$\Delta_{11}, \Delta_{12}, \Delta_{13}, \Delta_{14}.$$

Their definitions will be given below, when the corresponding terms are evaluated.

Define

$$\begin{aligned}
V_1 &= \mathbf{s}_{j_1}^* A_{j_1j_2}^{-1}(z)\mathbf{x}_n\mathbf{x}_n^* A_{j_1j_2}^{-1}(\bar{z})\mathbf{s}_{j_2},\\
V_3 &= \tilde{\mathbf{s}}_{j_2}^* \hat{A}_{j_1j_2}^{-1}(\bar{z})\hat{\mathbf{s}}_{j_2} - E(\hat{X}_{11}\tilde{X}_{11}^*)N^{-1}\operatorname{tr}\hat{A}_{j_1j_2}^{-1}(\bar{z})T_n,\\
V_2 &= \hat{\mathbf{s}}_{j_2}^* \hat{A}_{j_1j_2}^{-1}(\bar{z})\mathbf{x}_n\mathbf{x}_n^* \hat{A}_{j_1j_2}^{-1}(z)\hat{\mathbf{s}}_{j_1},\\
V_4 &= \hat{\mathbf{s}}_{j_1}^* \hat{A}_{j_1j_2}^{-1}(z)\tilde{\mathbf{s}}_{j_1} - E(\hat{X}_{11}\tilde{X}_{11}^*)N^{-1}\operatorname{tr}\hat{A}_{j_1j_2}^{-1}(z)T_n.
\end{aligned}$$

One can then easily prove that

$$\begin{aligned}
&E|V_3|^2 = o(N^{-1}), \qquad E|V_4|^2 = o(N^{-1}),\\
&E|V_1|^4 = O(N^{-4}) \quad\text{and}\quad E|V_2|^4 = O(N^{-4}).
\end{aligned}\tag{7.7}$$

Using (7.7) and (7.1), we obtain

$$\begin{aligned}
|\Delta_{11}| &= |E[\beta_{j_1}(z)\beta_{j_2}(\bar{z})\hat{\beta}_{j_1}(z)\hat{\beta}_{j_2}(\bar{z})\mathbf{s}_{j_1}^* A_{j_1j_2}^{-1}(z)\mathbf{x}_n\mathbf{x}_n^* A_{j_1j_2}^{-1}(\bar{z})\mathbf{s}_{j_2}\\
&\qquad \times \tilde{\mathbf{s}}_{j_2}^* \hat{A}_{j_1j_2}^{-1}(\bar{z})\hat{\mathbf{s}}_{j_2}\hat{\mathbf{s}}_{j_2}^* \hat{A}_{j_1j_2}^{-1}(\bar{z})\mathbf{x}_n\mathbf{x}_n^* \hat{A}_{j_1j_2}^{-1}(z)\hat{\mathbf{s}}_{j_1}\hat{\mathbf{s}}_{j_1}^* \hat{A}_{j_1j_2}^{-1}(z)\tilde{\mathbf{s}}_{j_1}]|\\
&\le K[E|V_1|^4]^{1/4}[E|V_2|^4]^{1/4}[E|V_3V_4|^2]^{1/2}\\
&\quad + KE|X_{11}|^2 I(|X_{11}| \ge \varepsilon_n n^{1/4})[E|V_1|^4]^{1/4}[E|V_2|^4]^{1/4}[E|V_3|^2]^{1/2}\\
&\quad + KE|X_{11}|^2 I(|X_{11}| \ge \varepsilon_n n^{1/4})[E|V_1|^4]^{1/4}[E|V_2|^4]^{1/4}[E|V_4|^2]^{1/2}\\
&\quad + K(E|X_{11}|^2 I(|X_{11}| \ge \varepsilon_n n^{1/4}))^2\\
&\qquad \times [E|\mathbf{s}_{j_1}^* A_{j_1j_2}^{-1}(z)\mathbf{x}_n\mathbf{x}_n^* \hat{A}_{j_1j_2}^{-1}(z)\hat{\mathbf{s}}_{j_1}|^2]^{1/2}\\
&\qquad \times [E|\hat{\mathbf{s}}_{j_2}^* \hat{A}_{j_1j_2}^{-1}(\bar{z})\mathbf{x}_n\mathbf{x}_n^* A_{j_1j_2}^{-1}(\bar{z})\mathbf{s}_{j_2}|^2]^{1/2}\\
&= o(N^{-3}),
\end{aligned}$$

where we make use of the independence between $\tilde{\mathbf{s}}_{j_2}^* \hat{A}_{j_1j_2}^{-1}(\bar{z})\hat{\mathbf{s}}_{j_2}$ and $\hat{\mathbf{s}}_{j_1}^* \hat{A}_{j_1j_2}^{-1}(z)\tilde{\mathbf{s}}_{j_1}$ when $\hat{A}_{j_1j_2}^{-1}(\bar{z})$ is given.

Let

$$V_5 = \hat{\mathbf{s}}_{j_2}^* \hat{A}_{j_1j_2}^{-1}(\bar{z})\mathbf{x}_n\mathbf{x}_n^* \hat{A}_{j_1j_2}^{-1}(z)\hat{\mathbf{s}}_{j_1}\mathbf{s}_{j_1}^* A_{j_1j_2}^{-1}(z)\mathbf{s}_{j_2}.$$

Similarly, we have

$$
\begin{aligned}
E|\mathbf{s}_{j_1}^* A_{j_1 j_2}^{-1}(z)\mathbf{s}_{j_2}^*|^4 &= O(N^{-2}),\\
E|V_5|^4 &= O(N^{-4}), \qquad (7.8)\\
E|\mathbf{s}_{j_2}^* A_{j_1 j_2}^{-1}(z)\mathbf{x}_n\mathbf{x}_n^* A_{j_1 j_2}^{-1}(\bar z)\mathbf{s}_{j_2}|^2 &= O(N^{-2}).
\end{aligned}
$$

Consequently, by (7.5), (7.7), (7.8) and (7.1), we have

$$
\begin{aligned}
|\Delta_{12}| &= |E[\beta_{j_1}(z)\beta_{j_2}(\bar z)\beta_{j_2 j_1}(z)\hat\beta_{j_1}(z)\hat\beta_{j_2}(\bar z)\mathbf{s}_{j_1}^* A_{j_1 j_2}^{-1}(z)\mathbf{s}_{j_2}\mathbf{s}_{j_2}^* A_{j_1 j_2}^{-1}(z)\mathbf{x}_n\\
&\qquad \times \mathbf{x}_n^* A_{j_1 j_2}^{-1}(\bar z)\mathbf{s}_{j_2}\tilde{\mathbf{s}}_{j_2}^*\hat A_{j_1 j_2}^{-1}(\bar z)\hat{\mathbf{s}}_{j_2}\\
&\qquad\qquad \times \hat{\mathbf{s}}_{j_2}^*\hat A_{j_1 j_2}^{-1}(\bar z)\mathbf{x}_n\mathbf{x}_n^*\hat A_{j_1 j_2}^{-1}(z)\hat{\mathbf{s}}_{j_1}\hat{\mathbf{s}}_{j_1}^*\hat A_{j_1 j_2}^{-1}(z)\tilde{\mathbf{s}}_{j_1}]|\\
&\le K[E|\mathbf{s}_{j_2}^* A_{j_1 j_2}^{-1}(z)\mathbf{x}_n\mathbf{x}_n^* A_{j_1 j_2}^{-1}(\bar z)\mathbf{s}_{j_2}|^2]^{1/2}\times[E|V_3V_4|^4]^{1/4}\times[E|V_5|^4]^{1/4}\\
&\quad + KE|X_{11}|^2 I(|X_{11}|\ge\varepsilon_n n^{1/4})[E|\mathbf{s}_{j_2}^* A_{j_1 j_2}^{-1}(z)\mathbf{x}_n\mathbf{x}_n^* A_{j_1 j_2}^{-1}(\bar z)\mathbf{s}_{j_2}|^2]^{1/2}\\
&\quad\quad \times[E|V_4|^4]^{1/4}[E|V_5|^4]^{1/4}\\
&\quad + KE|X_{11}|^2 I(|X_{11}|\ge\varepsilon_n n^{1/4})[E|\mathbf{s}_{j_2}^* A_{j_1 j_2}^{-1}(z)\mathbf{x}_n\mathbf{x}_n^* A_{j_1 j_2}^{-1}(\bar z)\mathbf{s}_{j_2}|^2]^{1/2}\\
&\quad\quad \times[E|V_3|^4]^{1/4}[E|V_5|^4]^{1/4}\\
&\quad + K(E|X_{11}|^2 I(|X_{11}|\ge\varepsilon_n n^{1/4}))^2[E|\mathbf{s}_{j_2}^* A_{j_1 j_2}^{-1}(z)\mathbf{x}_n\mathbf{x}_n^* A_{j_1 j_2}^{-1}(\bar z)\mathbf{s}_{j_2}|^2]^{1/2}\\
&\quad\quad \times[E|V_2|^4]^{1/4}[E|\mathbf{s}_{j_1}^* A_{j_1 j_2}^{-1}(z)\mathbf{s}_{j_2}^*|^4]^{1/4}\\
&= o(N^{-3}).
\end{aligned}
$$

Let $V_6 = s_{j_2}^* A_{j_1 j_2}^{-1}(z)\mathbf{x}_n\mathbf{x}_n^* A_{j_1 j_2}^{-1}(\bar z)\mathbf{s}_{j_1}\mathbf{s}_{j_1}^* A_{j_1 j_2}^{-1}(\bar z)\mathbf{s}_{j_2}$. By Lemma 1, one can obtain $E|V_6|^2 = O(N^{-3})$. Hence,

$$
\begin{aligned}
|\Delta_{13}| &= |E[\beta_{j_1}(z)\beta_{j_2}(\bar z)\beta_{j_1 j_2}(\bar z)\beta_{j_2 j_1}(z)\hat\beta_{j_1}(z)\hat\beta_{j_2}(\bar z)\mathbf{s}_{j_1}^* A_{j_1 j_2}^{-1}(z)\mathbf{s}_{j_2}\\
&\qquad \times \mathbf{s}_{j_2}^* A_{j_1 j_2}^{-1}(z)\mathbf{x}_n\mathbf{x}_n^* A_{j_1 j_2}^{-1}(\bar z)\mathbf{s}_{j_1}\mathbf{s}_{j_1}^* A_{j_1 j_2}^{-1}(\bar z)\mathbf{s}_{j_2}\tilde{\mathbf{s}}_{j_2}^*\hat A_{j_1 j_2}^{-1}(\bar z)\hat{\mathbf{s}}_{j_2}\\
&\qquad\qquad \times \hat{\mathbf{s}}_{j_2}^*\hat A_{j_1 j_2}^{-1}(\bar z)\mathbf{x}_n\mathbf{x}_n^*\hat A_{j_1 j_2}^{-1}(z)\hat{\mathbf{s}}_{j_1}\hat{\mathbf{s}}_{j_1}^*\hat A_{j_1 j_2}^{-1}(z)\tilde{\mathbf{s}}_{j_1}]|\\
&\le K[E|V_6|^2]^{1/2}[E|V_3V_4|^4]^{1/4}\times[E|V_5|^4]^{1/4}\\
&\quad + KE|X_{11}|^2 I(|X_{11}|\ge\varepsilon_n n^{1/4})K[E|V_6|^2]^{1/2}[E|V_4|^4]^{1/4}\times[E|V_5|^4]^{1/4}\\
&\quad + KE|X_{11}|^2 I(|X_{11}|\ge\varepsilon_n n^{1/4})[E|V_6|^2]^{1/2}[E|V_3|^4]^{1/4}\times[E|V_5|^4]^{1/4}\\
&\quad + K(E|X_{11}|^2 I(|X_{11}|\ge\varepsilon_n n^{1/4}))^2\\
&\quad\quad \times[E|V_6|^2]^{1/2}[E|V_2|^4]^{1/4}[E|\mathbf{s}_{j_1}^* A_{j_1 j_2}^{-1}(z)\mathbf{s}_{j_2}^*|^4]^{1/4}\\
&= o(N^{-3}).
\end{aligned}
$$

Before proceeding, we need the following estimate:

$$E|\mathbf{s}_p^* A_{pq}^{-1}(z)\mathbf{s}_q \hat{\mathbf{s}}_p^* \hat{A}_{pq}^{-1}(\bar{z})\mathbf{x}_n\mathbf{x}_n^* \hat{A}_{pq}^{-1}(z)\hat{\mathbf{s}}_q \hat{\mathbf{s}}_q^* \hat{A}_{pq}^{-1}(\bar{z})\hat{\mathbf{s}}_p|^4 = O(N^{-4}). \tag{7.9}$$

For $p \neq q$, write

$$\begin{aligned} T_n^{1/2} A_{pq}^{-1}(z) T_n^{1/2} &= (a_{uv}),\\ T_n^{1/2} \hat{A}_{pq}^{-1}(\bar{z})\mathbf{x}_n\mathbf{x}_n^* \hat{A}_{pq}^{-1}(z) T_n^{1/2} &= (b_{ml}),\\ T_n^{1/2} \hat{A}_{pq}^{-1}(\bar{z}) T_n^{1/2} &= (c_{ij}). \end{aligned}$$

Then

$$\begin{aligned} &E|\mathbf{s}_p^* A_{pq}^{-1}(z)\mathbf{s}_q \hat{\mathbf{s}}_p^* \hat{A}_{pq}^{-1}(\bar{z})\mathbf{x}_n\mathbf{x}_n^* \hat{A}_{pq}^{-1}(z)\hat{\mathbf{s}}_q \hat{\mathbf{s}}_q^* \hat{A}_{pq}^{-1}(\bar{z})\hat{\mathbf{s}}_p|^4\\ &\quad = \frac{1}{N^{12}} \sum E a_{u_1v_1} a_{u_2v_2} a_{u_3v_3} a_{u_4v_4}\\ &\qquad \times X_{u_1p}^* X_{u_2p} X_{u_3p}^* X_{u_4p} X_{v_1q} X_{v_2q}^* X_{v_3q} X_{v_4q}^*\\ &\qquad \times b_{m_1l_1} b_{m_2l_2} b_{m_3l_3} b_{m_4l_4} \qquad\qquad (7.10)\\ &\qquad \times \hat{X}_{m_1p}^* \hat{X}_{m_2p} \hat{X}_{m_3p}^* \hat{X}_{m_4p} \hat{X}_{l_1q} \hat{X}_{l_2q}^* \hat{X}_{l_3q} \hat{X}_{l_4q}^*\\ &\qquad \times c_{i_1j_1} c_{i_2j_2} c_{i_3j_3} c_{i_4j_4}\\ &\qquad \times \hat{X}_{i_1p} \hat{X}_{i_2p}^* \hat{X}_{i_3p} \hat{X}_{i_4p}^* \hat{X}_{j_1q}^* \hat{X}_{j_2q} \hat{X}_{j_3q}^* \hat{X}_{j_4q}, \end{aligned}$$

where the summations run over all indices except $p$ and $q$, and where $\hat{X}_{u_1p}^*$ denotes the complex conjugate and transpose of $\hat{X}_{u_1p}$, that is, the intercrossed element of the $u_1$th row and $p$th column of matrix $X_n$. Observe that:

(1) none of the indices of the above twenty-four random variables appear alone, otherwise the expectation is zero,

(2) $\operatorname{tr}(\hat{A}_j^{-1}(z)T_n\hat{A}_{pq}^{-1}(\bar{z})T_n) = O(N)$ and $|\hat{X}_{\boldsymbol{\cdot} p}| \le \varepsilon_n n^{1/4}$;

(3) $|\mathbf{x}_n^* \hat{A}_{pq}^{-1}(z)T_n\hat{A}_{pq}^{-1}(\bar{z})\mathbf{x}_n|$ is bounded. Combining the above, we obtain (7.9).

Similarly to (7.9), one can also verify that

$$\begin{aligned} &E|\hat{\mathbf{s}}_{j_2}^* \hat{A}_{j_1j_2}^{-1}(\bar{z})\hat{\mathbf{s}}_{j_1}\mathbf{s}_{j_2}^* A_{j_1j_2}^{-1}(z)\mathbf{x}_n\mathbf{x}_n^* A_{j_1j_2}^{-1}(\bar{z})\mathbf{s}_{j_1}\hat{\mathbf{s}}_{j_1}^* \hat{A}_{j_1j_2}^{-1}(\bar{z})\hat{\mathbf{s}}_{j_2}|^4 = O(N^{-4}),\\ &\qquad\qquad E|\tilde{\mathbf{s}}_{j_2}^* \hat{A}_{j_1j_2}^{-1}(\bar{z})\hat{\mathbf{s}}_{j_1}\hat{\mathbf{s}}_{j_2}^* \hat{A}_{j_1j_2}^{-1}(z)\tilde{\mathbf{s}}_{j_1}|^4 = o(N^{-2}) \end{aligned} \tag{7.11}$$

and that

$$E|\mathbf{s}_{j_1}^* A_{j_1j_2}^{-1}(\bar{z})\mathbf{s}_{j_2}\hat{\mathbf{s}}_{j_2}^* \hat{A}_{j_1j_2}^{-1}(\bar{z})\hat{\mathbf{s}}_{j_1}|^4 = O(N^{-2}). \tag{7.12}$$

It follows that

$$\begin{aligned}
|\Delta_{14}| &= |E[\beta_{j_1}(z)\beta_{j_2}(\bar z)\beta_{j_1j_2}(\bar z)\beta_{j_2j_1}(z)\hat\beta_{j_1}(z)\hat\beta_{j_2}(\bar z)(\hat\beta_{j_1j_2}(\bar z))^2(\hat\beta_{j_2j_1}(z))^2 \\
&\qquad \times \mathbf{s}_{j_1}^* A_{j_1j_2}^{-1}(z)\mathbf{s}_{j_2}\mathbf{s}_{j_2}^* A_{j_1j_2}^{-1}(z)\mathbf{x}_n\mathbf{x}_n^* A_{j_1j_2}^{-1}(\bar z)\mathbf{s}_{j_1}\mathbf{s}_{j_1}^* A_{j_1j_2}^{-1}(\bar z) \\
&\qquad \times \mathbf{s}_{j_2}\tilde{\mathbf{s}}_{j_2}^* \hat A_{j_1j_2}^{-1}(\bar z)\hat{\mathbf{s}}_{j_1}\hat{\mathbf{s}}_{j_1}^* \hat A_{j_1j_2}^{-1}(\bar z)\hat{\mathbf{s}}_{j_2}\hat{\mathbf{s}}_{j_2}^* \hat A_{j_1j_2}^{-1}(\bar z)\hat{\mathbf{s}}_{j_1}\hat{\mathbf{s}}_{j_1}^* \hat A_{j_1j_2}^{-1}(\bar z)\mathbf{x}_n \\
&\qquad\qquad \times \mathbf{x}_n^* \hat A_{j_1j_2}^{-1}(z)\hat{\mathbf{s}}_{j_2}\hat{\mathbf{s}}_{j_2}^* \hat A_{j_1j_2}^{-1}(z)\hat{\mathbf{s}}_{j_1}\hat{\mathbf{s}}_{j_1}^* \hat A_{j_1j_2}^{-1}(z)\hat{\mathbf{s}}_{j_2}\hat{\mathbf{s}}_{j_2}^* \hat A_{j_1j_2}^{-1}(z)\tilde{\mathbf{s}}_{j_1}]| \\
&\le K[E|\mathbf{s}_{j_1}^* A_{j_1j_2}^{-1}(z)\mathbf{s}_{j_2}\hat{\mathbf{s}}_{j_1}^* \hat A_{j_1j_2}^{-1}(\bar z)\mathbf{x}_n\mathbf{x}_n^* \hat A_{j_1j_2}^{-1}(z)\hat{\mathbf{s}}_{j_2}\hat{\mathbf{s}}_{j_2}^* \hat A_{j_1j_2}^{-1}(\bar z)\hat{\mathbf{s}}_{j_1}|^4]^{1/4} \\
&\quad \times [E|\hat{\mathbf{s}}_{j_1}^* \hat A_{j_1j_2}^{-1}(\bar z)\hat{\mathbf{s}}_{j_2}\mathbf{s}_{j_2}^* A_{j_1j_2}^{-1}(z)\mathbf{x}_n\mathbf{x}_n^* A_{j_1j_2}^{-1}(\bar z)\mathbf{s}_{j_1}\hat{\mathbf{s}}_{j_2}^* \hat A_{j_1j_2}^{-1}(z)\hat{\mathbf{s}}_{j_1}|^4]^{1/4} \\
&\quad \times [E|\hat{\mathbf{s}}_{j_1}^* \hat A_{j_1j_2}^{-1}(z)\hat{\mathbf{s}}_{j_2}\mathbf{s}_{j_1}^* A_{j_1j_2}^{-1}(\bar z)\mathbf{s}_{j_2}|^4]^{1/4} \\
&\quad \times [E|\tilde{\mathbf{s}}_{j_2}^* \hat A_{j_1j_2}^{-1}(\bar z)\hat{\mathbf{s}}_{j_1}\hat{\mathbf{s}}_{j_2}^* \hat A_{j_1j_2}^{-1}(z)\tilde{\mathbf{s}}_{j_1}|^4]^{1/4} \\
&= o(N^{-3}).
\end{aligned}$$

All other terms of $\Delta_1$ can be dealt with in a similar manner. Hence, combining the above arguments, one can conclude that

$$n\sum_{j_1\ne j_2}\Delta_1 \to 0.$$

Similarly, one can also prove that $n\sum_{j_1\ne j_2}\Delta_2$, $n\sum_{j_1\ne j_2}\Delta_3$ and $n\sum_{j_1\ne j_2}\Delta_4$ converge to zero. Hence, $\omega_{12}\to 0$. Combining this with (7.6), we have

$$\omega_1 \xrightarrow{i.p.} 0.$$

Using the same method, one can prove that $\omega_2$ converges to zero in probability. Hence, we complete the truncation and centralization of $X_{ij}$. The rescaling of random variables $X_{ij}$ can be completed in a similar way. By Theorem 1 in [2], we have

$$\liminf_{n\to\infty}\min(u_r-\lambda_{\max}^{\hat A_n}, \lambda_{\min}^{\hat A_n}-u_l)>0 \qquad \text{a.s.}$$

and

$$\liminf_{n\to\infty}\min(u_r-\lambda_{\max}^{A_n}, \lambda_{\min}^{A_n}-u_l)>0 \qquad \text{a.s.}$$

Hence, for $z\in\mathcal{C}_l, u_l>0$ or $z\in\mathcal{C}_r$, we can also proceed with the truncation, centralization and rescaling of $X_{ij}$ by the above method. Thus the proof is complete. □

**Acknowledgments.** The authors would like to thank the referee for valuable comments. Dr. G. M. Pan wishes to thank the members of the department for their hospitality when he visited the Department of Statistics and Applied Probability at National University of Singapore.

Z. D. Bai
Key Laboratory for Applied Statistics of MOE
and College of Mathematics and Statistics
Northeast Normal University
Changchun, 130024
China
E-mail: baizd@nenu.edu.cn
stabaizd@nus.edu.sg

B. Q. Miao
Department of Statistics and Finance
University of Science
and Technology of China
Hefei, 230026
China
E-mail: bqmiao@ustc.edu.cn

G. M. Pan
Department of Statistics and Finance
University of Science
and Technology of China
Hefei, 230026
China
E-mail: stapgm@gmail.com